\documentclass[a4paper]{siamart171218}
\usepackage[percent]{overpic}
\usepackage[numbers]{natbib}
\usepackage{amsmath}
\usepackage{amscd}
\usepackage{amssymb}
\usepackage{graphicx}
\usepackage{url}
\usepackage{color}
\usepackage{listings}
\usepackage{booktabs}
\usepackage{multirow}
\usepackage{textcomp}
\usepackage{mathtools}
\usepackage{subcaption}
\usepackage{hyperref}
\usepackage[utf8]{inputenc}
\usepackage{tikz}
\usetikzlibrary{trees,calc}
\usepackage{pgfplots}
\usetikzlibrary{plotmarks}
\usepackage{pgfplotstable}
\usepackage{subcaption}
\usepackage{environ}
\usepackage{algorithm}
\usepackage{enumerate}
\usepackage[noend]{algpseudocode}
\usepackage{color}
\usepackage{xspace}

\definecolor{blue}{rgb}{0.2980392156862745, 0.4470588235294118, 0.6901960784313725}
\definecolor{green}{rgb}{0.3333333333333333, 0.6588235294117647, 0.40784313725490196}
\definecolor{red}{rgb}{0.7686274509803922, 0.3058823529411765, 0.3215686274509804}

\renewcommand{\Re}{\ensuremath{\mathrm{Re}}\xspace}
\newcommand{\honev}{\ensuremath{{H}^1(\Omega; \mathbb{R}^d)}\xspace}
\newcommand{\ltwov}{\ensuremath{{L}^2(\Omega; \mathbb{R}^d)}\xspace}
\newcommand{\ltwo}{\ensuremath{{L}^2(\Omega)}\xspace}

\newcommand{\eps}[1]{\ensuremath{\varepsilon(#1)}}

\newcommand{\ds}{\ \mathrm{d}s}
\newcommand{\dx}{\ \mathrm{d}x}
\newcommand{\Pq}{\ensuremath{\mathrm{P}_{Q_h}}}

\newcommand{\Ptwo}{\ensuremath{\mathbb{P}_2}\xspace}
\newcommand{\Pthree}{\ensuremath{\mathbb{P}_3}\xspace}
\newcommand{\PtwoPzero}{\ensuremath{[\mathbb{P}_2]^2\mathrm{-}\mathbb{P}_0}\xspace}
\newcommand{\PthreePzero}{\ensuremath{[\mathbb{P}_3]^3\mathrm{-}\mathbb{P}_0}\xspace}
\newcommand{\Pzero}{\ensuremath{\mathbb{P}_0}\xspace}
\newcommand{\Pv}{\ensuremath{\mathbb{P}_v}\xspace}

\newcommand{\PoneFB}{\ensuremath{\mathbb{P}_1 \oplus B^F_3}\xspace}
\newcommand{\PtwoFB}{\ensuremath{\mathbb{P}_2 \oplus B^F_3}\xspace}

\newcommand{\advect}[2]{\ensuremath{(#2 \cdot \nabla) #1}}
\newcommand{\mesh}{\ensuremath{\mathcal{M}}\xspace}

\newtheorem{remark}[theorem]{Remark}

\definecolor{darkblue}{rgb}{0.00,0.00,0.55}
\definecolor{black}{rgb}{0.00,0.00,0.00}
\hypersetup{
    linkcolor = darkblue,
    anchorcolor = darkblue,
    citecolor = darkblue,
    filecolor = darkblue,
    urlcolor = darkblue,
    colorlinks  = true,
}

\title{An augmented Lagrangian preconditioner for the 3D stationary incompressible Navier--Stokes equations at high Reynolds number
\thanks{Submitted to the editors \today.
\funding{This research is supported by the Engineering and Physical Sciences
Research Council [grant numbers EP/K030930/1, EP/M011151/1, EP/M011054/1,
and EP/L000407/1], and by the EPSRC Centre For Doctoral Training in
Industrially Focused Mathematical Modelling [grant number EP/L015803/1] in
collaboration with London Computational Solutions. LM also acknowledges
support from the UK Fluids Network [EPSRC grant number EP/N032861/1] for
funding a visit to Oxford. This work used the ARCHER UK National
Supercomputing Service (\url{http://www.archer.ac.uk}). The authors would
like to acknowledge useful discussions with M.~Benzi, M.~A.~Olshanskii, A.~J.~Wathen, D.~J.~Silvester and
J.~Sch\"oberl, to thank M.~A.~Olshanskii for supplying the code for the
Oseen solver described in \cite{benzi2006}, and to thank M.~G.~Knepley for
assistance with PETSc.}}}

\author{
  Patrick~E.~Farrell\thanks{Mathematical Institute, University of Oxford, Oxford, UK (\email{patrick.farrell@maths.ox.ac.uk}).}
  \and
  Lawrence Mitchell\thanks{Department of Computer Science, Durham University, Durham, UK
    (\email{lawrence.mitchell@durham.ac.uk}).}
  \and
  Florian~Wechsung\thanks{Mathematical Institute, University of Oxford,
    Oxford, UK (\email{florian.wechsung@maths.ox.ac.uk}). }} 

\headers{Preconditioners for high-Re 3D stationary flow}{P.~E.~Farrell, L. Mitchell, and F.~Wechsung}

\begin{document}
\maketitle

\begin{abstract}
In Benzi \& Olshanskii (SIAM J.~Sci.~Comput., 28(6) (2006)) a preconditioner of
augmented Lagrangian type was presented for the two-dimensional stationary
incompressible Navier--Stokes equations that exhibits convergence almost
independent of Reynolds number. The algorithm relies on a highly specialized multigrid
method involving a custom prolongation operator and for robustness requires the
use of piecewise constant finite elements for the pressure. However, the
prolongation operator and velocity element used do not directly extend to three
dimensions: the local solves necessary in the prolongation operator do not
satisfy the inf-sup condition. In this work we generalize the preconditioner to
three dimensions, proposing alternative finite elements for the velocity and
prolongation operators for which the preconditioner works robustly.  The solver
is effective at high Reynolds number: on a three-dimensional lid-driven
cavity problem with approximately one billion degrees of freedom, the average
number of Krylov iterations per Newton step varies from 4.5 at $\Re = 10$ to 3
at $\Re = 1000$ and 5 at $\Re = 5000$.
\end{abstract}



\section{Introduction} \label{sec:introduction}


We consider the stationary incompressible Newtonian Navier--Stokes equations:
given a bounded Lipschitz domain $\Omega \subset \mathbb{R}^d$, $d \in \{2,
3\}$, find $(u, p) \in \honev \times Q$ such that
\begin{subequations} \label{eqn:ns}
\begin{alignat}{2}
    -  \nabla\cdot2\nu\eps{u} + \advect{u}{u} + \nabla p &= f \quad && \text{ in } \Omega, \label{eqn:momentum} \\
\nabla \cdot u &= 0 \quad && \text{ in } \Omega, \\
             u &= g \quad && \text{ on } \Gamma_D, \\
            2\nu \eps{u} \cdot n &= pn \quad && \text{ on } \Gamma_N,
\end{alignat}
\end{subequations}
where $\eps{u} = \frac{1}{2}(\nabla u + \nabla u^T)$, $\nu >0$ is the kinematic viscosity, $f \in \ltwov$, $n$ is the outward-facing unit normal to $\partial
\Omega$, $\Gamma_D$ and $\Gamma_N$ are disjoint with $\Gamma_D \cup \Gamma_N =
\partial \Omega$, and $g \in H^{1/2}(\Gamma_D; \mathbb{R}^d)$. If $|\Gamma_N| > 0$, then
a suitable trial space for the pressure is $Q := \ltwo$; if $|\Gamma_N| = 0$,
then the pressure is only defined up to an additive constant and $Q :=
L^2_0(\Omega)$ is used instead. The Reynolds number, defined as $\Re = \frac{UL}{\nu}$ where $U$ is the characteristic
velocity and $L$ is the characteristic length scale of the flow, is an important dimensionless number governing
the nature of the flow.
The Navier--Stokes
equations are of enormous practical importance in science and
industry, but are very difficult to solve, especially for
large Reynolds number.
The importance of these equations has motivated a great deal of
research on algorithms for their solution; for a general overview of the field,
see the textbooks of Turek \cite{turek1999}, Elman, Silvester \& Wathen
\cite{elman2014}, or Brandt \cite{brandt2011}.

%
%
%
%

After Newton linearization and a suitable spatial discretization of
\eqref{eqn:ns},
nonsymmetric linear systems of saddle point type must be solved:
\begin{equation} \label{eqn:sp}
\begin{pmatrix}
A & B^T \\
B & 0
\end{pmatrix}
\begin{pmatrix}
\delta u \\ \delta p
\end{pmatrix}
=
\begin{pmatrix}
b \\ c
\end{pmatrix},
\end{equation}
where $A$ is the discrete linearized momentum operator, $B^T$ is the discrete gradient
operator, $B$ is the discrete divergence operator, and $\delta u$ and $\delta p$
are the updates to the coefficients for velocity and pressure respectively. One
strategy to solve these systems is to employ a monolithic multigrid iteration on
the entire system with a suitable coupled relaxation method, such as the algorithms of
Vanka \cite{vanka1986} or Braess \& Sarazin \cite{braess1997}. Vanka iteration
works well for moderate Reynolds numbers \cite{turek1999}, but the iteration
counts have been observed to increase significantly once the Reynolds number
becomes large \cite{benzi2006}.

An alternative approach to solving \eqref{eqn:sp} is to build
preconditioners based on block factorizations
\cite{murphy2000,ipsen2001,benzi2005,elman2014,wathen2015}. This strategy can be grounded
in an insightful
functional analytic framework that guides the development of solvers whose
convergence is independent of parameter values and mesh size $h$, at least in the case where
\eqref{eqn:sp} is symmetric \cite{mardal2011}. Block Gaussian elimination reduces the problem of solving the coupled
linear system to that of solving smaller separate linear systems involving the
matrix $A$ and the Schur complement $S = -B A^{-1} B^T$. If a fast solver is
available for $A$, the main difficulty is solving linear systems involving $S$,
as this matrix is generally dense and cannot be stored explicitly for large
problems. Tractable approximations $\tilde{S}^{-1}$ to $S^{-1}$ must be devised
on a PDE-specific basis.

For the Stokes equations, the Schur complement is spectrally equivalent to the
viscosity-weighted pressure mass matrix \cite{silvester1994}. For the
Navier--Stokes equations this choice yields mesh-independent convergence and is
effective for very small Reynolds numbers, but the convergence deteriorates badly with Reynolds
number \cite{elman1996,elman2014}. The pressure convection-diffusion (PCD)
approach \cite{kay2002} constructs an auxiliary convection-diffusion operator on
the pressure space, and hypothesizes that a certain commutator is small. This
yields an approximation to the Schur complement inverse that involves the
inverse of the Laplacian on the pressure space, the application of the
auxiliary convection-diffusion operator, and the inverse of the pressure mass matrix. The
least-squares commutator (LSC) approach \cite{elman2006} is based on a similar
idea, but derives the commutator algebraically. Both of these
approaches perform well for moderate Reynolds numbers. Numerical experiments
comparing the performance of our approach to these algorithms are provided
in section \ref{sec:numerical}.

In 2006, Benzi \& Olshanskii proposed an augmented Lagrangian approach for
controlling the Schur complement of \eqref{eqn:sp}
\cite{olshanskii2003,benzi2006,olshanskii2008,benzi2011b,fortin1983,kobelkov1995}. The
idea, referred to as grad-div stabilization, is to introduce an additional term in
the equations that does not change the continuous solution, but does modify the
Schur complement. The continuous form of the stabilization replaces
\eqref{eqn:momentum} with
\begin{alignat}{2} \label{eqn:nsal}
- \nabla\cdot2\nu\eps{u} + \advect{u}{u} + \nabla p  - \gamma \nabla \nabla \cdot u &= f \quad && \text{ in } \Omega,
\end{alignat}
for $\gamma > 0$.
As $\nabla \cdot u = 0$, the solutions of \eqref{eqn:nsal} and
\eqref{eqn:momentum} are the same. The discrete variant of this approach
replaces \eqref{eqn:sp} with
\begin{equation} \label{eqn:spal}
\begin{pmatrix}
A + \gamma B^T M_p^{-1} B & B^T \\
B & 0
\end{pmatrix}
\begin{pmatrix}
\delta u \\ \delta p
\end{pmatrix}
=
\begin{pmatrix}
b + \gamma B^T M_p^{-1} c \\ c
\end{pmatrix},
\end{equation}
where $M_p$ is the pressure mass matrix. This modified system has the same discrete solutions as \eqref{eqn:sp}, as $B\delta u = 0$.
With either variant, for $\gamma$ not too small, the Schur
complement inverse is well approximated by
\begin{equation} \label{eqn:schur}
S^{-1} \approx -(\nu + \gamma) M_p^{-1},
\end{equation}
where $M_p$ is the pressure mass matrix. This approximation improves
as $\gamma$ increases (section \ref{sec:augmentedlagrangian}). In either case, we denote the discretized augmented
momentum block as $A_\gamma$.

\begin{remark}
  The continuous form of the grad-div stabilization has some further
  appealing characteristics. For example, it significantly improves
  the pressure-robustness of discretizations where the
  incompressibility constraint is not enforced pointwise
  \cite{olshanskii2003,heister2012,john2017}. It also arises in other
  contexts in the numerical analysis of \eqref{eqn:ns}. For example,
  Boffi \& Lovadina \cite{boffi1997} showed that the addition of the
  term $h^{-1/2} (\nabla \cdot u, \nabla \cdot v)_{\ltwo}$ to the weak
  form of the \PtwoPzero discretization of \eqref{eqn:ns} improves its
  convergence order. It also arises in the iterated penalty
  \cite{temam1968,brenner2008} and artificial compressibility
  \cite{chorin1967} methods for the Stokes and Navier--Stokes
  equations.
\end{remark}

The tradeoff with either variant of this approach is that developing fast
solvers for $A_\gamma$ becomes significantly more difficult.
The divergence operator has a large kernel (the range of the curl operator) and
hence standard multigrid relaxation methods are ineffective.

A key insight of Benzi \& Olshanskii was that a specialized multigrid algorithm
could be built for $A_\gamma$ \cite{benzi2006,olshanskii2008} by applying the seminal work of Sch\"oberl
\cite{schoberl1999}. The algorithm combines four
ingredients, each of which is crucial to the effectiveness of the method: (i)
the discrete variant of the grad-div stabilization; (ii) streamline-upwind
Petrov-Galerkin (SUPG) stabilization of the advective term; (iii) a multigrid
relaxation that effectively treats errors in the kernel of the discrete
divergence term; (iv) a specialized prolongation operator whose continuity
constant is independent of $\gamma$ and $\nu$. This scheme exhibits outer
iteration counts that grow only very slowly with Reynolds number
\cite{benzi2006}. However, it is described as difficult to implement
\cite{hamilton2010,benzi2011}, and so most of the works that use
grad-div stabilization and the Schur complement approximation \eqref{eqn:schur}
employ either matrix factorization as the inner solver
\cite{deniet2007,urrehman2008,borm2010,he2011,heister2012} or a block-triangular
approximation to $A_\gamma$ \cite{benzi2011,hamilton2010,benzi2011b, he2018}. This block-triangular
approximation decouples linear systems involving $A_\gamma$ into $d$ scalar anisotropic
advection--diffusion problems, which may be solved with algebraic multigrid
techniques. However, this simplicity comes at a price; the scheme is much more
sensitive to the choice of $\gamma$, and its convergence deteriorates somewhat
as the Reynolds number increases \cite{benzi2011}.

The main contribution of this paper is the extension of the robust multigrid scheme for
the inner velocity problem arising in the augmented Lagrangian preconditioner to
three dimensions. The previous work of Benzi \& Olshanskii only considered the
case $d=2$. While the same general strategy applies in three dimensions, the extension
is nontrivial: if the finite elements used in \cite{schoberl1999,benzi2006}
are applied in three dimensions, the prolongation operator involves the
solution of ill-posed local problems. We propose appropriate finite element
discretizations and matching prolongation operators that exhibit Reynolds-robust
iteration counts in three dimensions.

A second contribution is the release of an open-source parallel implementation of the
solver in two and three dimensions, built on Firedrake \cite{rathgeber2016} and PETSc \cite{balay2017}.
This has required substantial modifications to Firedrake, PETSc,
UFL \cite{alnaes2012} and TSFC \cite{homolya2017}, as well as minor
developments in FIAT \cite{kirby2004}. The solver heavily relies on and
extends the solver infrastructure developed in \cite{kirby2018},
enabling easy composition and nesting of preconditioners in PETSc and
Firedrake. To express the local solves involved in the
relaxation and prolongation operator, we have developed a new
preconditioner in PETSc that allows for the simple expression of
general additive subspace correction methods. For example, the same code
that does patchwise relaxation can be used to formulate line
smoothers, plane smoothers, or Vanka iteration, and we expect that it
will be of substantial interest for other applications as well.

The remainder of the paper is laid out as follows. The discretization
and the grad-div stabilization are described in section
\ref{sec:discretization}. The augmented Lagrangian approach is
explained in section \ref{sec:augmentedlagrangian}. The multigrid
cycle for the augmented momentum block is described in section
\ref{sec:solvingmomentum}. Numerical experiments analyzing its
performance and comparing it to PCD and LSC are reported in section
\ref{sec:numerical}. Finally, conclusions and prospects for future
improvements are given in section \ref{sec:conclusions}.

\section{Formulation and discretization} \label{sec:discretization}

For boundary data $g \in H^{1/2}(\Gamma_D)$, let
\begin{equation}
V_g = \{v \in H^1(\Omega; \mathbb{R}^d) : \left. v \right|_{\Gamma_D} = g \}.
\end{equation}
The initial weak form of \eqref{eqn:ns} is: find $(u, p) \in V_g \times Q$ such that
\begin{equation}
    \int_\Omega 2\nu\eps{u}: \nabla v \dx + \int_\Omega \advect{u}{u} \cdot v \dx - \int_\Omega p \nabla \cdot v \dx - \int_\Omega q \nabla \cdot u \dx = \int_\Omega f \cdot v \dx,
\end{equation}
for all $(v, q) \in V_0 \times Q$.
This will be extended before discretization in two ways. The first is a consistent SUPG stabilization; it is well
known that straightforward Galerkin discretizations of advection-dominated problems are
oscillatory \cite{brooks1982,turek1999,quarteroni2008,elman2014}. In addition,
it is widely observed that mesh-dependent SUPG stabilization is highly advantageous for multigrid
smoothers on advection-dominated problems \cite{ramage1999,turek1999}. The strong form
of the momentum residual is given by
\begin{equation}
\mathcal{L}(u, p) = - \nabla \cdot 2\nu\eps{u} + \advect{u}{u} + \nabla p - f
\end{equation}
and the following term is added to the weak form:
\begin{equation}
\int_\Omega \delta(u) \mathcal{L}(u, p) \cdot \big(\advect{v}{u}\big) \dx.
\end{equation}
Here $\delta(u)$ is a weighting function that should be small in regions where
the flow is well-resolved and large where stabilization is necessary.
The particular form employed in this work is
\begin{equation}
\delta(u) = \delta_d \left( \frac{4 \|u\|^2}{h^2} + \frac{144\nu^2}{h^4} \right)^{-1/2},
\end{equation}
with $\delta_d = 1$ in two dimensions and $\delta_d = 1/20$ in three dimensions.
To the best of our knowledge this form was first suggested in
\cite[eq.~(3.58)]{shakib1991}. It is important to take account of the dependence of
$\delta$ on the (unknown) solution $u$ when taking the derivatives required by
Newton's method; in this work, these derivatives are calculated automatically
and symbolically by the Unified Form Language \cite{alnaes2012}.

The second modification is the augmented Lagrangian term described
above in \eqref{eqn:nsal} and \eqref{eqn:spal}. If the
continuous variant is employed, the term
\begin{equation} \label{eqn:ctsgraddiv}
\gamma \int_\Omega \nabla \cdot u \, \nabla \cdot v \dx
\end{equation}
is added to the weak form, while if the discrete variant is employed, the term
\begin{equation} \label{eqn:disgraddiv}
\gamma \int_\Omega \big(\Pq \nabla \cdot u\big) \nabla \cdot v \dx
\end{equation}
is added instead, where $\Pq: L^2(\Omega) \to Q_h$ is the projection operator onto the discrete
pressure space $Q_h$.
The continuous grad-div stabilization changes the discrete solution
computed if the discrete velocity $u_h$ does not satisfy $\nabla \cdot
u_h = 0$ pointwise,
whereas the discrete variant does not.
The effect of the continuous approach is to penalize $\|\nabla \cdot
u_h\|_{L^2(\Omega)}$ and thereby improve the discrete enforcement of the
incompressibility constraint \cite{olshanskii2003,heister2012,john2017}.
Nevertheless, in this work we use the discrete variant \eqref{eqn:disgraddiv}.
The reason for this is that the kernel of $\Pq \mathrm{div}$ is much more
straightforward to characterize than the kernel of $\mathrm{div}$ if
$Q_h$ is chosen to be the space of piecewise constants:
\begin{equation}
Q_h(\mesh) = \{q \in L^2(\Omega) : \left.q\right|_K \in \mathbb{P}_0(K) \ \forall K \ \in \mesh\},
\end{equation}
where $\mesh$ is a simplicial mesh of the domain $\Omega$.
By the divergence theorem,
$u_h \in \mathrm{ker}(\Pq
\mathrm{div})$ if and only if for any $K \in \mesh$, $u_h$ satisfies
\begin{equation}
\int_{\partial K} \! u_h \cdot n \ds = 0.
\end{equation}
This characterization will be extremely useful for dealing with errors in the
kernel in the multigrid relaxation, as it ensures that the kernel is spanned
by basis functions with local support \cite[\S VI.8]{brezzi1991}. Note also that this choice of $Q_h$ removes
the dependency of $\mathcal{L}$ on $p$ (as $\nabla p$ is zero on each element),
thereby eliminating any extra contribution to the top-right block of the linearized
system to be solved, thus preserving the symmetry between the top-right and bottom-left
blocks of the matrix.

After these modifications, the final discrete weak form to be solved is:
find $(u, p) \in \big(V_h \cap V_g\big) \times \big(Q_h \cap
Q\big)$ such that
\begin{align}
& \ \int_\Omega 2\nu\eps{u} : \nabla v \dx + \int_\Omega \advect{u}{u} \cdot v\dx \nonumber
+ \int_\Omega \delta(u) \mathcal{L}(u) \cdot \big(\advect{v}{u}\big) \dx \nonumber\\
& + \gamma \int_\Omega \big(\Pq \nabla \cdot u\big) \nabla \cdot v \dx
 - \int_\Omega p \nabla \cdot v \dx - \int_\Omega q \nabla \cdot u \dx =
   \int_\Omega f \cdot v \dx, \label{eqn:finalweak}
\end{align}
for all $(v, q) \in \big(V_h \cap V_0\big) \times \big(Q_h \cap Q\big)$, with
the choice of $V_h$ to be discussed below.
\subsection{Choice of velocity space}
\label{sec:velocityspace}
We now turn our attention to choosing an appropriate space for the discrete
velocities. Define the space $V_h(\mesh)$ used for the velocity as
\begin{equation}
V_h(\mesh) = \{ v \in H^1(\Omega) : \left.v\right|_K \in \Pv(K) \ \forall \ K \in \mesh \},
\end{equation}
for some choice of $\Pv(K)$. The first condition on $\Pv(K)$ is that $V_h$ must
be inf-sup stable when combined with $Q_h$ for the pressure. Unfortunately,
both in two and in three dimensions, the $[\mathbb{P}_1]^d-\mathbb{P}_0$ element
combining piecewise linear functions for the velocity space together with
piecewise constants for the pressure does not satisfy the inf-sup condition on general meshes.
This means that the velocity space needs to be enriched, with the resulting
element pairs (e.g.~$\PtwoPzero$ and $\PthreePzero$) exhibiting a suboptimal
convergence rate.

In \cite{benzi2006} the element pair
$[\mathbb{P}_1\text{iso}\,\mathbb{P}_2]^2\mathrm{-}\mathbb{P}_0$ is
used, which is obtained by considering a $[\mathbb{P}_1]^2$ element on
a once refined mesh for the velocity. This element has the same
number of degrees of freedom as $[\Ptwo]^2\mathrm{-}\Pzero$.
Neither $[\Ptwo]^3\mathrm{-}\Pzero$ nor $[\mathbb{P}_1\text{iso}\,\Ptwo]^3\mathrm{-}\Pzero$
are inf-sup stable on a single regularly refined tetrahedron,
which as we shall see in
section~\ref{sec:prolongation} is crucial for the effectiveness of the
preconditioner. They are missing degrees of freedom on the facets of the
tetrahedra which are necessary to stabilize the jump of the pressure
field.

\begin{figure}[htbp]
\centering
\includegraphics[width=0.24\linewidth]{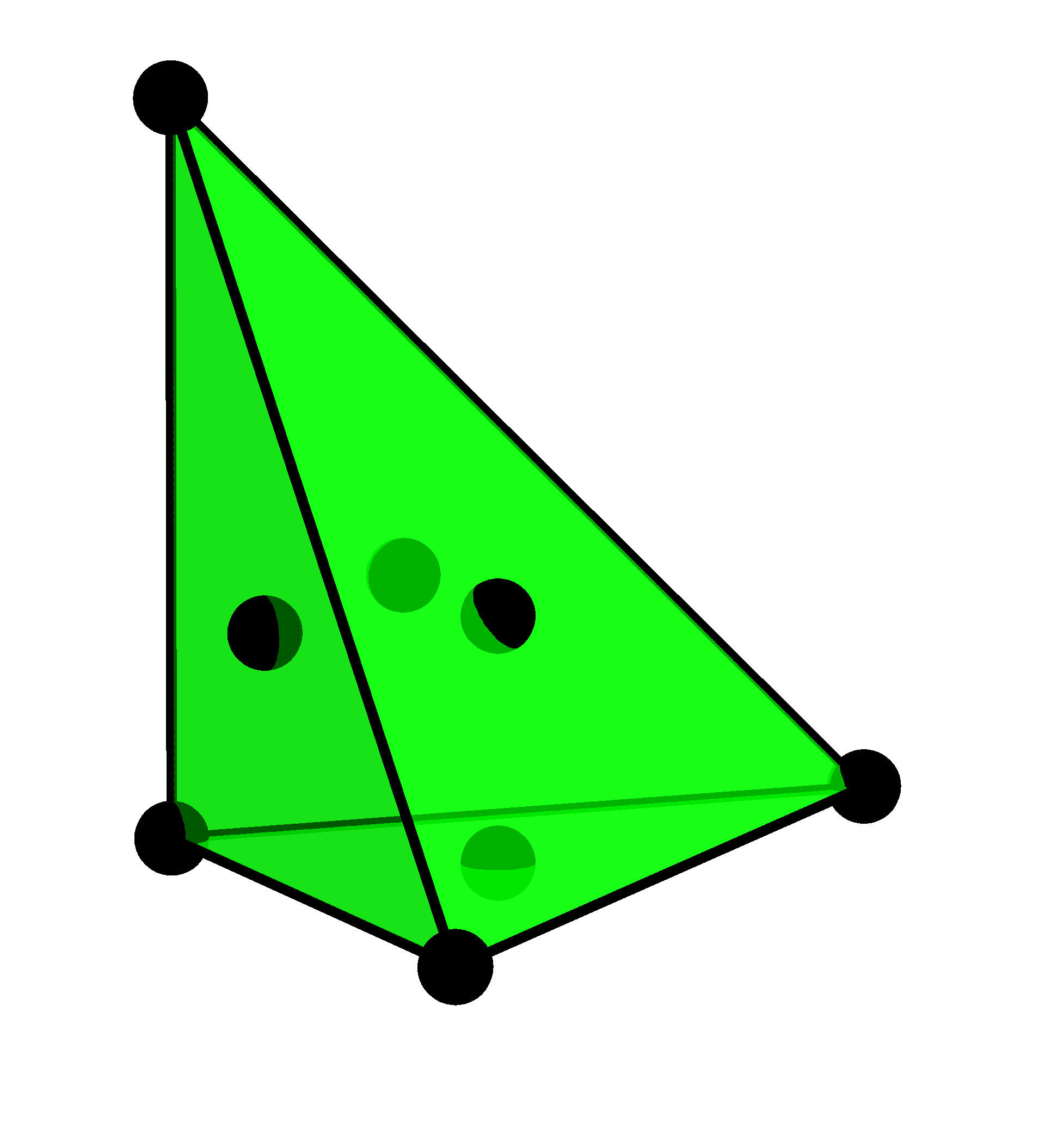}
\includegraphics[width=0.24\linewidth]{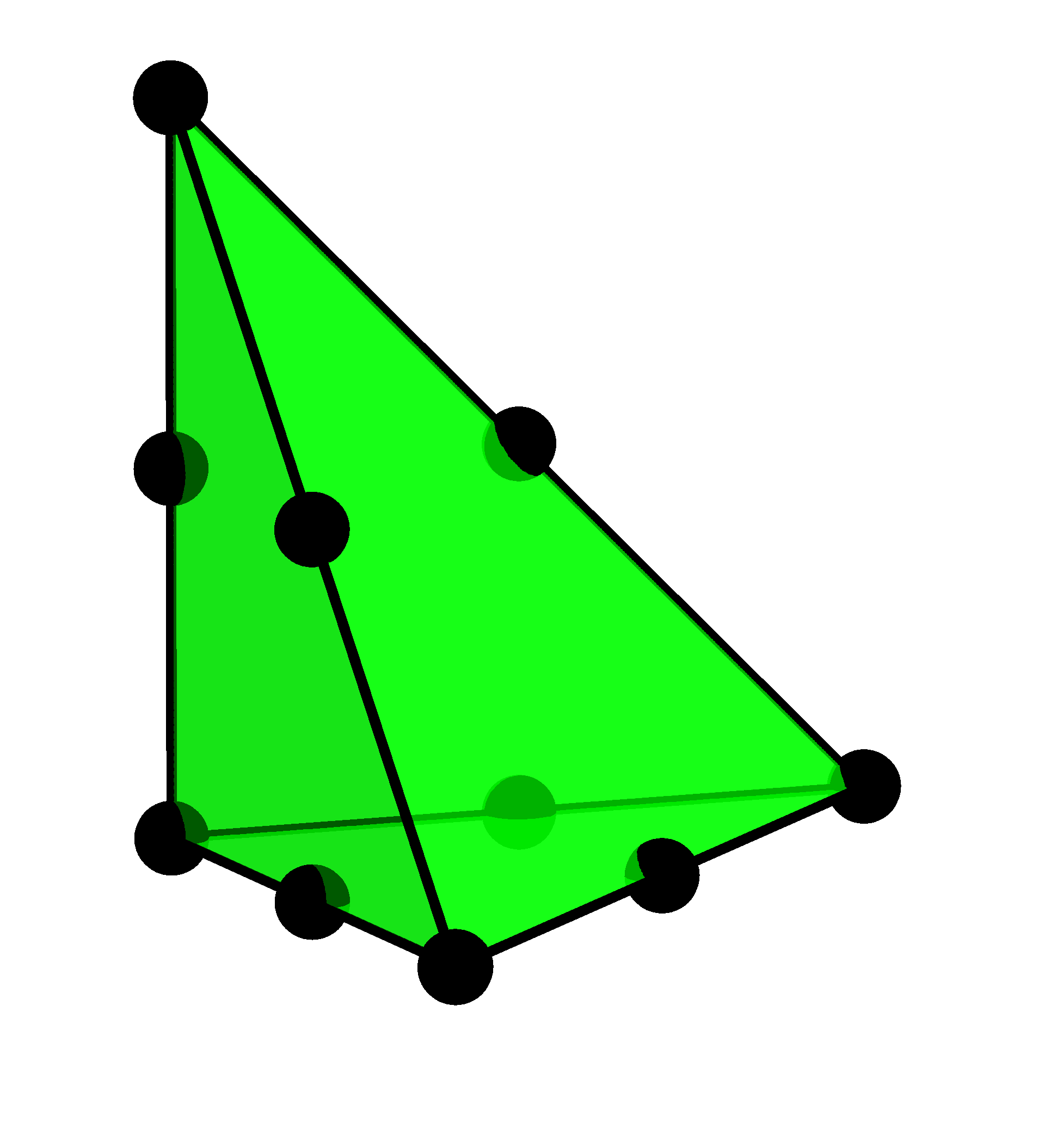}
\includegraphics[width=0.24\linewidth]{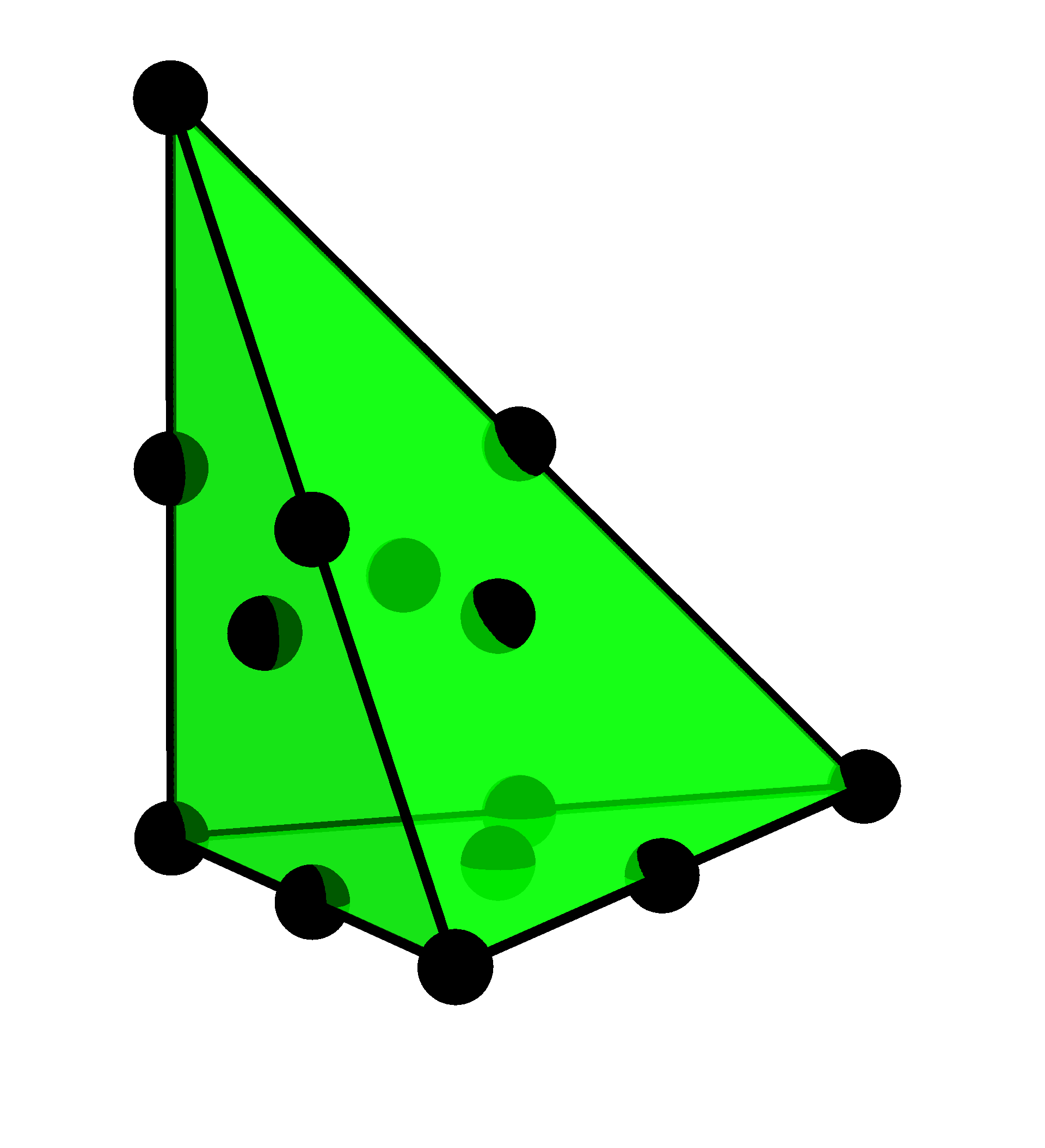}
\includegraphics[width=0.24\linewidth]{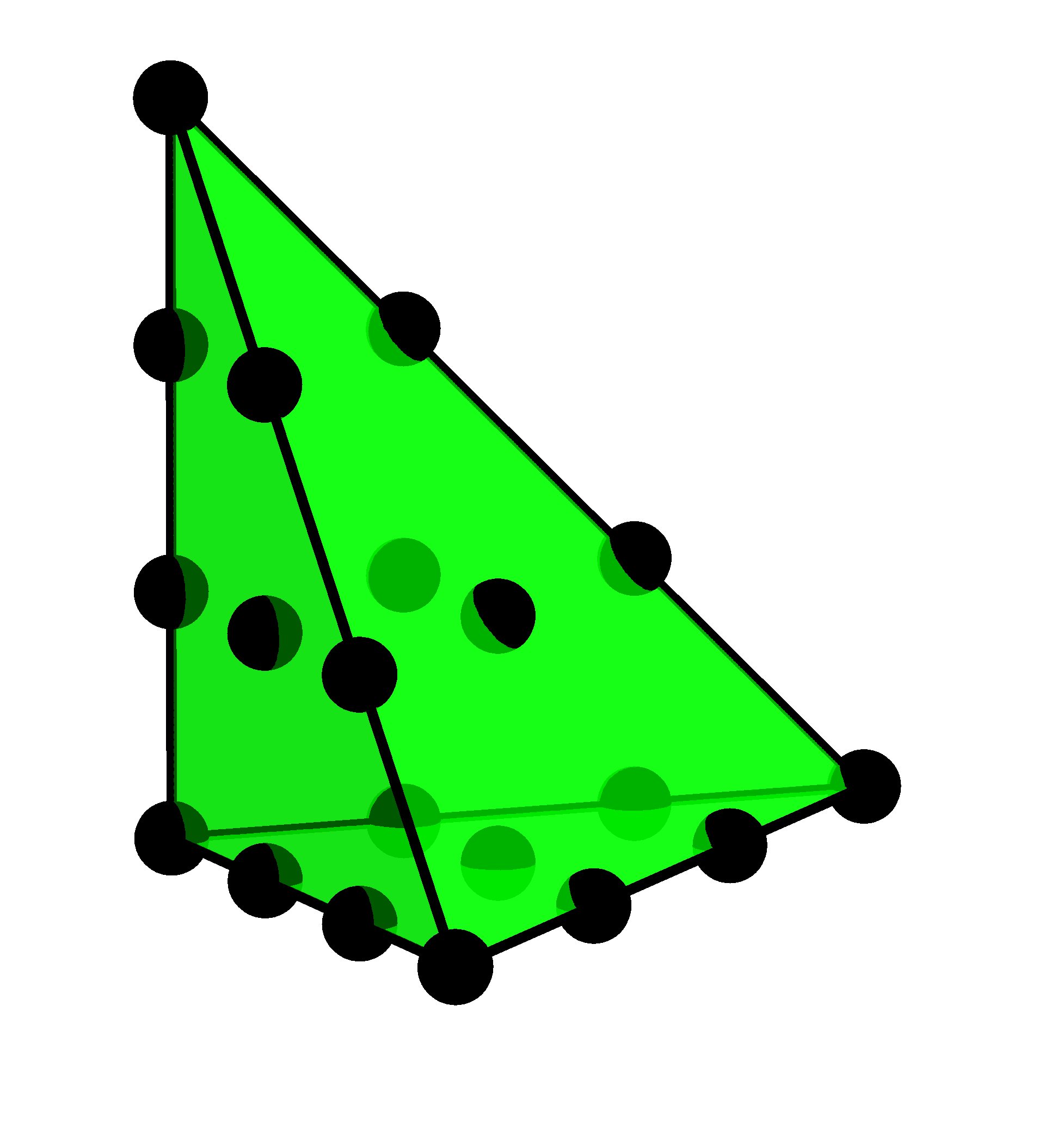}
\caption{The \PoneFB, \Ptwo, \PtwoFB and \Pthree elements.}
\label{fig:elements}
\end{figure}

Increasing the degree of the velocity space to piecewise cubic
polynomials, i.e.~choosing the element pair
$[\Pthree]^3\mathrm{-}\mathbb{P}_0$, introduces additional degrees of
freedom on the facets and results in a stable element pair. However, this
element is extremely expensive while being suboptimal by two orders for the
velocity. Alternatively, Bernardi \& Raugel \cite{bernardi1985, bernardi1985b} suggest enriching the
piecewise linear velocity space with bubble functions on each facet\footnote{The bubble
function on each facet is the product of the barycentric coordinates that are nonzero on that facet.}. While it
is only necessary to add a single bubble function for the normal
component of the velocity on each facet, this adds significant complexity to the implementation
as these functions are not affine equivalent; they require a Piola transform
to preserve the normal orientation. This means that the basis functions
associated with vertices and those associated with facets need to be
pulled-back differently, complicating the implementation. For this reason
we choose instead to enrich the space with facet bubbles for all three components of the velocity,
obtaining the $[\PoneFB]^3\mathrm{-}\Pzero$ element. As can be seen in
Figure~\ref{fig:elements}, this results in an element with significantly fewer
degrees of freedom than $[\Pthree]^3\mathrm{-}\Pzero$. We also show the $[\PtwoFB]^3$ element in
Figure~\ref{fig:elements}; we will demonstrate in
Section~\ref{sec:prolongation} that these elements satisfy a particular property
that is useful in the prolongation.

\begin{remark}
Pressure elements other than $\Pzero$ have been considered for the augmented
Lagrangian preconditioner.
Benzi \& Olshanskii \cite[Table 6.2]{benzi2006} also present results for the
$[\mathbb{P}_1\text{iso}\,\mathbb{P}_2]^2\mathrm{-}\mathbb{P}_1$ pair,
where the pressure mass matrix solve in $\Pq$ is approximated by the inverse of a
diagonal matrix. However, for this element pairing the developed multigrid scheme is not independent of the ratio $\gamma/\nu$ and hence as $\nu$ decreases, $\gamma$ has to be decreased correspondingly.
This in turn leads to worse control of the Schur complement and consequent growth in iteration counts.
\end{remark}
\begin{remark}
Effective smoothers for partial differential equations involving the continuous
term $(\nabla \cdot u, \nabla \cdot v)_{L^2(\Omega)}$ have been proposed in
other contexts \cite{arnold1997,hiptmair1997}; it should be possible to extend these
approaches to the two- and three-dimensional advection-dominated case
considered here, and thereby enable the use of finite elements with
advantageous properties such as optimal convergence rates and exact enforcement
of the incompressibility constraint.
\end{remark}
\section{The augmented Lagrangian method} \label{sec:augmentedlagrangian}


The Schur complement of the matrix in \eqref{eqn:spal} is given by
\begin{equation}
    S = -B(A+\gamma B^T M_p^{-1} B)^{-1} B^T.
\end{equation}
From the Sherman–Morrison–Woodbury formula it follows (e.g.~\cite[Theorem 3.2]{bacuta2006}) that
\begin{equation}
  S^{-1} = -(BA^{-1} B^T)^{-1} - \gamma M_p^{-1}.
\end{equation}
From this we obtain immediately that
\begin{equation}
    \tilde S^{-1} = - (\nu + \gamma) M_p^{-1}
\end{equation}
is a good approximation to $S^{-1}$ as $\gamma\to\infty$ \emph{for fixed mesh
size and viscosity}. To understand the quality of the approximation for
finite $\gamma$ as $\nu$ or $h$ change, we need to consider how $\nu M_p^{-1}$
approximates $(BA^{-1}B^T)^{-1}$.

It is well known that the eigenvalues of a matrix do not characterize the
convergence of GMRES for a linear system \cite{greenbaum1996}. Instead, it is
necessary to bound the field-of-values of the preconditioned system
\cite[Theorem 3.2]{starke1997}, \cite[Corollary 6.2]{eiermann2001}, \cite[\S 1.3]{olshanskii2014}. This
analysis was performed by Benzi \& Olshanskii \cite{benzi2011b} for both the
ideal and the modified augmented Lagrangian preconditioner for the Oseen
problem, using general results of Loghin \& Wathen \cite{loghin2004}. One of
the key ingredients in this analysis is that the momentum operator is coercive
with constant $\nu$. They use this to prove that the choice of
$\gamma \sim \nu^{-1}$ results in an optimal preconditioner (assuming exact
solves of the momentum block). However, it is well known
\cite[p.~300]{girault1986} that the momentum operator of the
Newton linearization of \eqref{eqn:ns} is only coercive for $\nu > \nu_0$ for
some problem-dependent $\nu_0$. Fortin \& Glowinski remark
\cite[p.~85]{fortin1983} that this is typically a very restrictive condition:
for $\nu > \nu_0$ the Stokes approximation itself is adequate. This proof
strategy would therefore require significant extension to apply to the Newton
linearization considered here.

In practice, Benzi \& Olshanskii \cite{benzi2006} observe that a constant choice of $\gamma$ yields
mesh-independent and essentially Reynolds number independent results. As our
multigrid solver for the momentum block is robust with respect to $\gamma$, we simply choose
$\gamma$ large. In the experiments of section \ref{sec:numerical}, we take the value $\gamma = 10^4$,
to match the largest Reynolds number considered.

\section{Solving the augmented momentum block} \label{sec:solvingmomentum}
The key challenge with the augmented Lagrangian strategy is the solution of the
augmented momentum block $A_\gamma$. The grad-div term has a large nullspace,
rendering standard relaxation methods (point-block Jacobi or Gauss--Seidel)
ineffective as $\nu \to 0$. In this
section we explain the specialized multigrid algorithm of Benzi \& Olshanskii,
along with the modifications required to extend the method to three dimensions.
The multigrid method has two components: a $\nu$- and $\gamma$-robust relaxation
method, and a kernel-preserving prolongation operator. In subsections
\ref{sec:relaxation} and \ref{sec:prolongation} we first consider the augmented
Stokes momentum operator without the linearized advection terms, to study in the
simplest possible situation the difficulties arising with the grad-div term. We
then comment on the case with advection in subsection \ref{sec:advection}.

To understand the properties required of the relaxation and
prolongation, it suffices to consider a two-level scheme. We use the
subscripts $h$ and $H$ to denote function spaces, bilinear forms, and
meshes on the fine and coarse levels respectively.

\subsection{Relaxation} \label{sec:relaxation}
The augmented Stokes momentum problem is of the form: find $u \in V_{h,0} := V_h
\cap V_0$ such that
\begin{equation} \label{eqn:bilinearstokes}
a_h(u, v) \coloneqq (2\nu\eps{u}, \nabla v) + \gamma (\Pq \nabla \cdot u, \nabla \cdot
v) = (f, v)
\end{equation}
for all $v \in V_{h,0}$.
The viscosity term is symmetric and positive definite; the discrete grad-div term
is positive semidefinite.
As $\nu \to 0$ or $\gamma \to \infty$ this system becomes nearly
singular and standard relaxation methods such as Gauss--Seidel or Jacobi
iterations perform poorly. The essential difficulty is in computing the
component of the solution in the kernel
\begin{equation}
\mathcal{N}_h := \{u \in V_{h,0} : (\Pq \nabla \cdot u, \nabla \cdot v) = 0
\ \forall\ v \in V_{h,0}\}
\end{equation}
of the grad-div term.
Sch\"oberl \cite[Theorem~4.1]{schoberl1999b} and Lee et
al.~\cite[Theorem~4.2]{lee2007} consider subspace
correction methods for this class of problem.
The key result of these works is that if a subspace decomposition
\begin{equation}
V_{h,0} = \sum_i V_i
\end{equation}
satisfies the kernel decomposition property
\begin{equation} \label{eqn:kernel-decomposition}
\mathcal{N}_h = \sum_i \left(V_i \cap \mathcal{N}_h \right)
\end{equation}
then the resulting subspace correction method (a block Gauss--Seidel
or Jacobi iteration) is robust with respect to $\nu$ and
$\gamma$. This is why characterizing the kernel of the
grad-div term is crucial, and why the discrete variant is easier
to solve: the kernel $\mathcal{N}_h$ is spanned by basis functions with
local support around each vertex.

More specifically, for each vertex $v_i$ in the mesh \mesh, its
\emph{star} is the patch of elements sharing $v_i$:
\begin{equation}
K_i \coloneqq \bigcup_{{K \in \mesh \,: \,v_i \in K}} K.
\end{equation}
The subspace decomposition is then given by
\begin{equation}
V_i \coloneqq \{ \phi_j \in V_{h,0} : \mathrm{supp}(\phi_j) \subset K_i \}.
\end{equation}
We call the resulting patchwise block relaxation method a \emph{star iteration}.
Note that homogeneous Dirichlet conditions are imposed on the boundary of each
star patch.
This relaxation method has been employed for robust multigrid methods in $H(\mathrm{div})$
and $H(\mathrm{curl})$ \cite{arnold2000}.

For the reader's convenience, we briefly summarize the argument of \cite[Section~4.1.2]{schoberl1999b} to see why this decomposition satisfies \eqref{eqn:kernel-decomposition}.
Observe that a discretely divergence-free vector field $u_h\in \mathcal{N}_h$
can be suitably modified in the \emph{interior} of each cell to become
continuously divergence-free by solving a local Stokes problem.
Denote this continuously divergence-free vector field by $\tilde u$ and recall
that then $\tilde u = \nabla \times \phi$ for some vector field $\phi$.
Choosing a partition of unity $\{\rho_i\}_i$ with $\sum_i \rho_i=1$ and
$\mathrm{supp}(\rho_i) \subset K_i$ we define $\phi_i = \rho_i \phi$ and obtain
a decomposition
\begin{equation}
    \phi = \sum_i \phi_i.
\end{equation}
For such a partition of unity to exist, every point in the mesh has to be in the
interior of at least one patch. The decomposition of the mesh into
star patches is the smallest decomposition of a finite element mesh that
satisfies this property.
Now let $\Pi_1:V\to V_h$ be a Scott--Zhang interpolation operator using facet averaging \cite{scott1990}.
Then it holds that $\Pi_1(u_h) = u_h$.
Furthermore, define $\Pi_2:V \to V_h$ as in the classical proof for inf-sup stability of the $\PtwoPzero$ element \cite[Proposition~3.1]{boffi2008}:
\begin{equation}
  \begin{aligned}
    \Pi_2(v)(M) = 0, \quad \text{ for all vertices } M,\\
    \int_{F}\Pi_2(v) \ds = \int_F v \ds, \quad \text{ for all facets } F.
  \end{aligned}
\end{equation}
Now define $I(v) = \Pi_1(v) + \Pi_2(v-\Pi_1(v))$, then it holds that
\begin{equation}
  \begin{aligned}
    I(v_h) &= v_h \quad \text{ for all } v_h\in V_h,\\
    \int_F I(v) \ds &= \int_F v \ds \quad \text{ for all } v\in V.\\
  \end{aligned}
\end{equation}
Now define $u_i=I(\nabla \times \phi_i)$ and conclude that
\begin{equation}
  \sum_i u_i = \sum_i I (\nabla \times \phi_i) = I(\nabla \times \phi) = I(\tilde u)= u_h.
\end{equation}
Lastly, using the fact that we are considering piecewise constant pressures, $u_i \in V_i\cap \mathcal{N}_h$ follows from
\begin{equation}
  \int_{\partial K} u_i \ds = \int_{\partial K} \nabla \times \phi_i \ds = \int_K \nabla\cdot(\nabla\times\phi_i) \dx = 0.
\end{equation}
\subsection{Prolongation} \label{sec:prolongation}
The second key ingredient of the multigrid scheme is the prolongation operator
that maps $V_H$ to $V_h$. To get an intuition for the properties required, let
$E_H: V_H \to V_h$ be the prolongation operator obtained by interpolating a
function $u_H\in V_H$ at the degrees of freedom of $V_h$.
The continuity of $E_H$ in the energy norm induced by the bilinear form $a_h$
defined in \eqref{eqn:bilinearstokes} is a key assumption in Sch\"oberl's proof
of the optimality of a two-level multigrid scheme \cite[Lemma
3.5]{schoberl1999b}. In order for the scheme to be robust, this continuity
constant must be uniform in $\nu$ and $\gamma$. Calculating, we observe that
\begin{equation}
    \begin{aligned}
        \|u_H \|_{a_H}^2 &= \nu\| \eps{u_H}\|_{L^2}^2 + \gamma \|P_{Q_H} (\nabla \cdot u_H)\|_{L^2}^2\\
        \|E_H u_H \|_{a_h}^2 &= \nu \| \eps{E_H u_H}\|_{L^2}^2 + \gamma \|P_{Q_h} (\nabla \cdot (E_Hu_H))\|_{L^2}^2.
    \end{aligned}
\end{equation}
The key difficulty lies in the second term of this norm. To see this, observe
that for an element $u_H \in \mathcal{N}_H$ the second term in
$\|u_H\|_{a_H}^2$ vanishes, but since it does not necessarily hold that $E_H u_H \in  \mathcal{N}_h$,
the corresponding term in $\|E_Hu_H\|_{a_h}^2$ might be
large.

To avoid this, we must modify the prolongation operator to map fields that are
discretely divergence-free on the coarse grid to fields that are (nearly)
discretely divergence-free on the fine grid.

To begin, we assume that there is a decomposition $Q_h = \tilde Q_H \oplus Q_T$ and a subspace $V_T\subset V_h$ that satisfies $V_T \subset \ker(P_{Q_H}(\nabla\cdot))$.
Sch\"oberl proved that if the pairing $V_T \mathrm{-} Q_T$ satisfies the inf-sup condition and if
\begin{alignat}{2}
    (P_{Q_h}(\nabla\cdot(E_H u_H)), \tilde q_H)_{L^2} &= (P_{Q_H} (\nabla\cdot u_H), \tilde q_H)_{L^2} && \quad\forall u_H\in V_H, \tilde q_H \in \tilde Q_H\label{eqn:decomposition-exterior}\\
    (P_{Q_h}(\nabla\cdot u_T), \tilde q_H)_{L^2} &= 0 &&\quad \forall u_T\in V_T,\ \tilde q_H \in \tilde Q_H\label{eqn:decomposition-interior}
\end{alignat}
then the prolongation $\tilde E_H$ defined as
\begin{equation} \label{eqn:schoberlprolong}
    \tilde E_H u_H = E_H u_h - w_T,
\end{equation}
where $w_T\in V_T$ satisfies
\begin{equation}\label{eqn:local-prolongation-problem}
    a_h(w_T, v_T) = a_h(E_H u_H, v_T) \quad \forall v_T \in V_T,
\end{equation}
is continuous in the energy norm. The continuity constant is uniform in $\nu$ and $\gamma$.
In this case, the decomposition of $Q_h$ is chosen as
\begin{alignat}{2}
    \tilde Q_H &\coloneqq Q_H \\
    Q_T &\coloneqq \{ q_h \in Q_h : P_{Q_H}(q_h) = 0\}
\end{alignat}
and we choose
\begin{equation}
    V_T \coloneqq \{ v_h \in V_h : v_h\vert_{\partial K} = 0 \ \forall K \in \mesh_H\}.
\end{equation}
The idea behind this is the following: \eqref{eqn:decomposition-exterior}
guarantees that prolongation preserves the flux across \emph{coarse grid}
facets. Then a correction term $w_T\in V_T$ that corrects the flux across the fine grid
facets is subtracted. The condition \eqref{eqn:decomposition-interior} guarantees that this
correction does not affect the flux across the coarse facets.

\begin{remark} The
definition of $V_T$ implies that the problem in
\eqref{eqn:local-prolongation-problem} can be solved locally on each coarse
grid element. This is crucial for an efficient implementation.
\end{remark}

\begin{remark}
Decompositions $\tilde Q_H \neq Q_H$ arise in other problems, such as
in Reissner--Mindlin plates \cite[Section 4.2.2]{schoberl1999b}.
\end{remark}

\begin{figure}[htbp]
    \centering
    \includegraphics[width=0.5\textwidth]{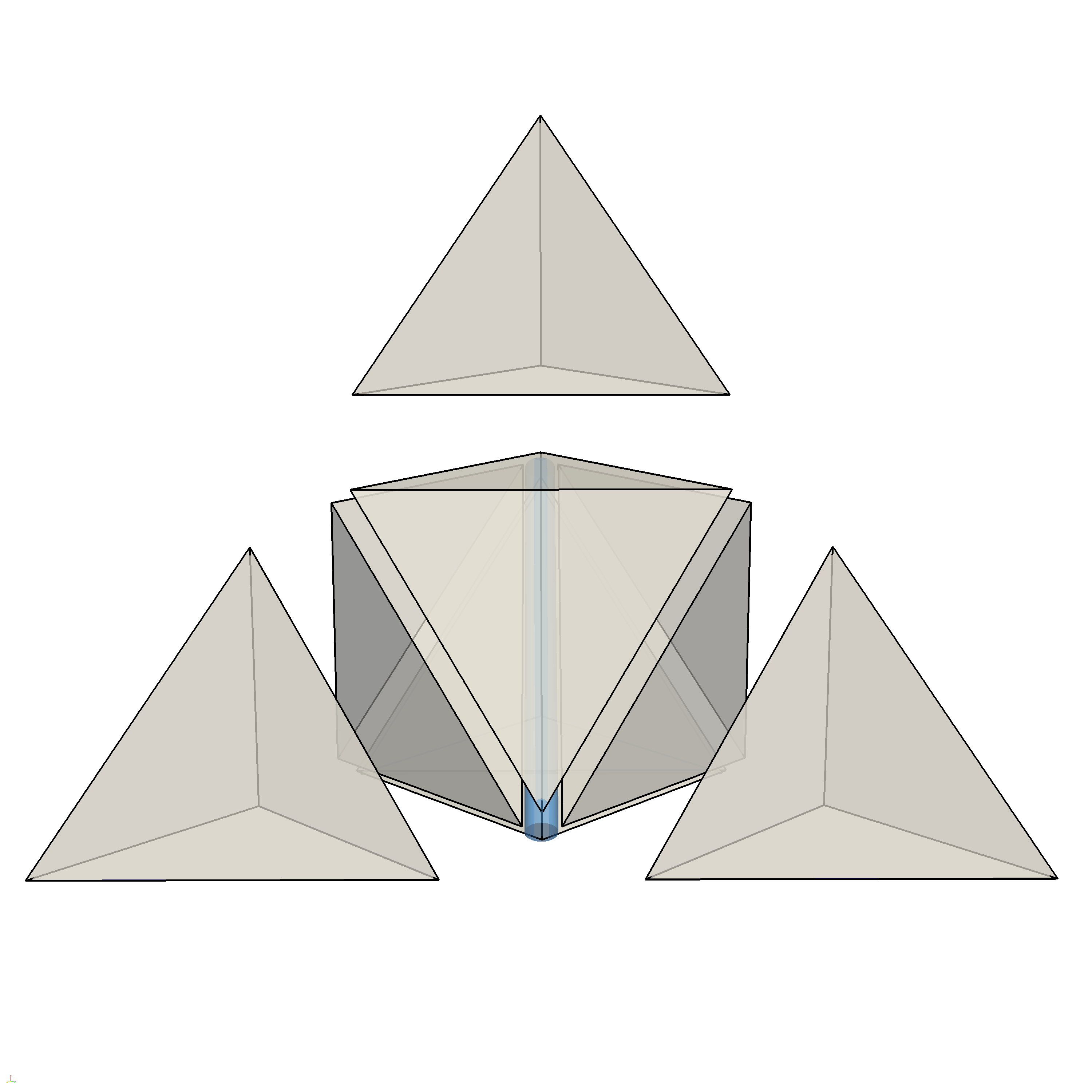}
    \caption{The uniform refinement of a coarse mesh tetrahedron yields eight fine mesh tetrahedra. Only the edge highlighted in blue does not lie
    entirely on the boundary of the tetrahedron.}
    \label{fig:exploded-tet}
\end{figure}

In \cite{schoberl1999, benzi2006} the $[\Ptwo]^2-\Pzero$ element is used. For
this element choice it holds that $V_H\subset V_h$ and hence
\eqref{eqn:decomposition-interior} is satisfied trivially. However, in
three dimensions the pairing $V_T\mathrm{-}Q_T$ resulting from the
choice $[\Ptwo]^3-\Pzero$ is not inf-sup stable. This can easily seen by
counting degrees of freedom: $[\Ptwo]^3$ only has degrees of freedom on edges and
vertices. Since there are zero vertices and only one edge not entirely on the boundary
of the refined coarse tetrahedron (see Figure~\ref{fig:exploded-tet}), we have $\dim(V_T)=3$. On the other hand, the
pressure space satisfies $\dim(Q_T) = 7$ (one dimension is fixed by the
nullspace). The local solve can therefore not be well-posed.

The choice $[\Pthree]^3-\Pzero$ alleviates this problem of ill-posedness on the
coarse cell and still satisfies $V_H\subset V_h$. However, as described in
Section~\ref{sec:velocityspace}, this element is quite expensive without
improving accuracy of the solution.

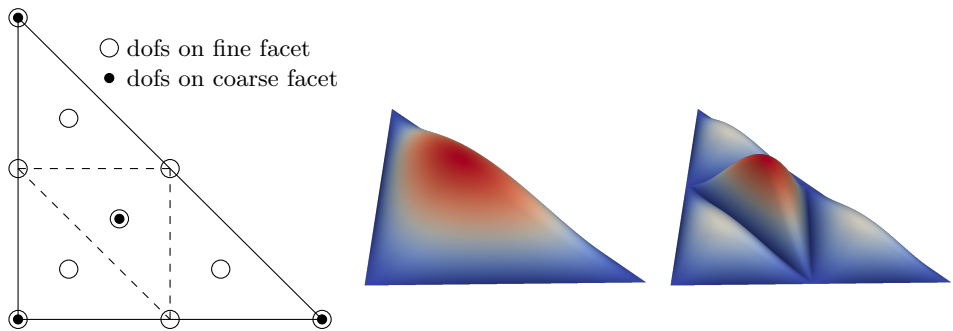
\begin{figure}[htbp]
    \centering
    \begin{tikzpicture}[scale=4.0]
        \draw (0, 0) -- (1, 0);
        \draw (1, 0) -- (0, 1);
        \draw (0, 1) -- (0, 0);
        \draw[dashed] (0.5, 0) -- (0.5, 0.5);
        \draw[dashed] (0.5, 0.5) -- (0.0, 0.5);
        \draw[dashed] (0.0, 0.5) -- (0.5, 0.0);
        \draw[fill=black]  (1./3., 1./3.) circle (0.015);
        \draw[fill=black]  (0, 0) circle (0.015);
        \draw[fill=black]  (0, 1) circle (0.015);
        \draw[fill=black]  (1, 0) circle (0.015);
        \draw  (1./3., 1./3.) circle (0.03);
        \draw  (2./3., 1./6.) circle (0.03);
        \draw  (1./6., 2./3.) circle (0.03);
        \draw  (1./6., 1./6.) circle (0.03);
        \draw  (1, 0) circle (0.03);
        \draw  (0, 0) circle (0.03);
        \draw  (0, 1) circle (0.03);
        \draw  (0.5, 0) circle (0.03);
        \draw  (0.5, 0.5) circle (0.03);
        \draw  (0, 0.5) circle (0.03);

        \draw  (0.3, 0.9) circle (0.03);
        \draw[anchor=west] node at (0.3, 0.9) {\small \ dofs on fine facet};
        \draw[fill=black]  (0.3, 0.8) circle (0.015);
        \draw[anchor=west] node at (0.3, 0.8) {\small \ dofs on coarse facet};
    \end{tikzpicture}
    \includegraphics[width=0.3\textwidth]{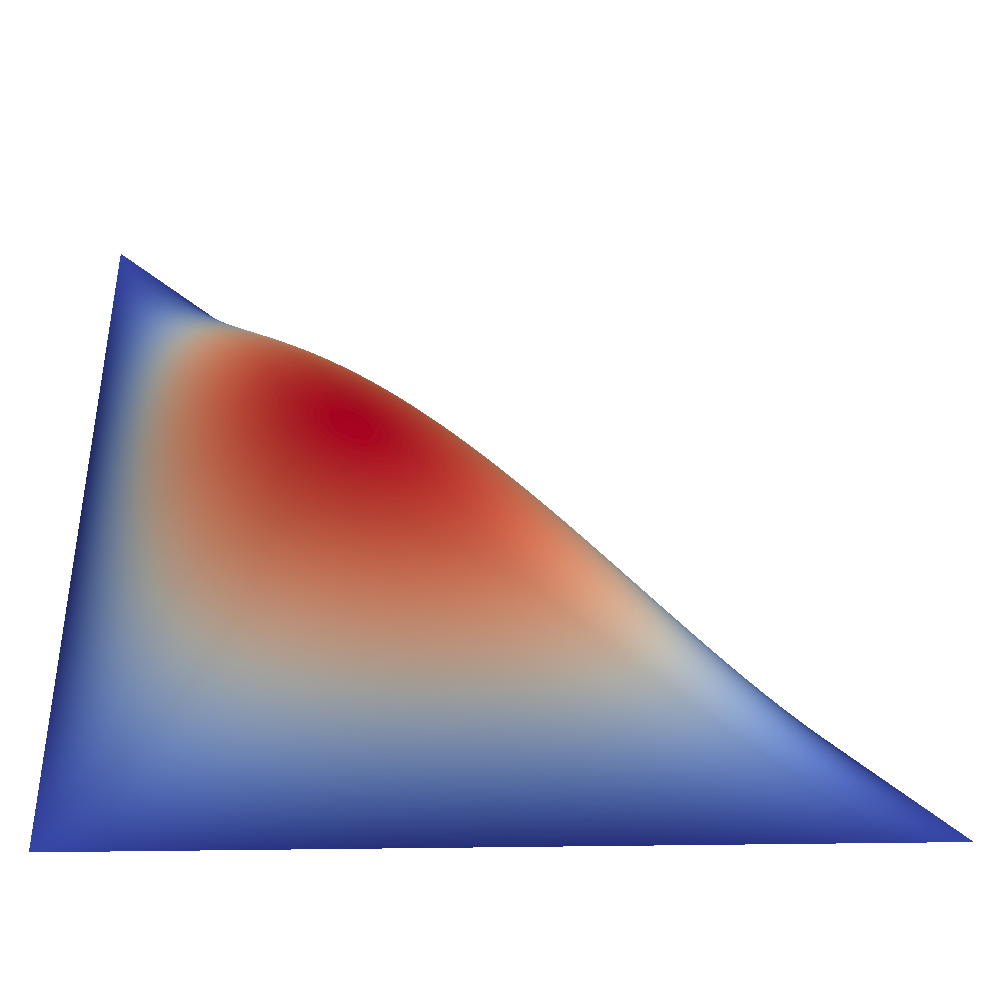}
    \includegraphics[width=0.3\textwidth]{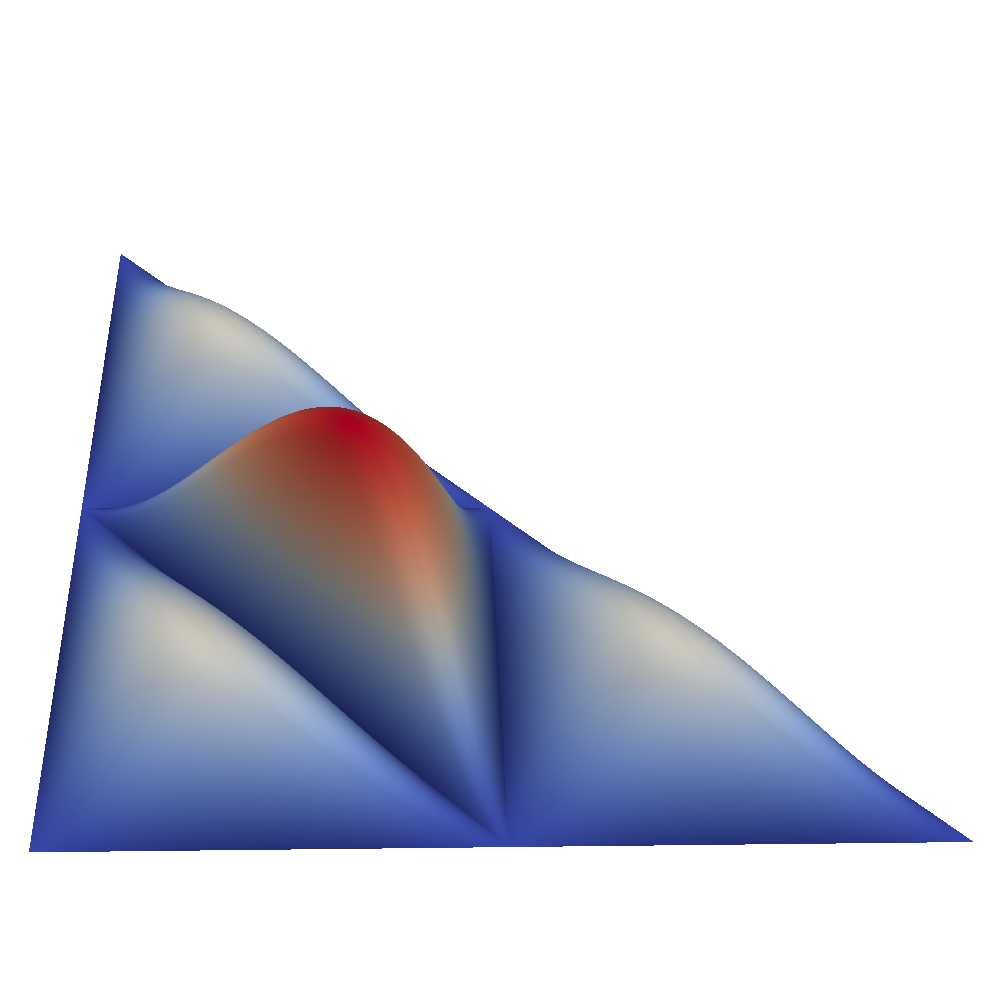}
    \caption{Left: Degrees of freedom on the facet of a coarse cell and its refinement. Middle: Bubble function on a coarse facet. Right: Prolongation of a bubble function.}
    \label{fig:fb-prolongation}
\end{figure}

A much cheaper alternative
is offered by the $[\PoneFB]^3\mathrm{-}\Pzero$ element. This does satisfy
the inf-sup condition but
violates $V_H\subset V_h$. The
non-nestedness is demonstrated in Figure~\ref{fig:fb-prolongation}; a coarse
bubble cannot be interpolated exactly by functions in $V_h$. In particular,
the flux across the coarse grid faces is not preserved, hence violating
\eqref{eqn:decomposition-exterior}.
A brief calculation shows that every coarse
grid bubble is interpolated by four fine grid bubbles: one with
coefficient $1$, the other three with coefficient $1/2$. From this it follows
immediately that the integral of the prolonged bubble is equal to $(1 + 3
\cdot \frac{1}2)/4 = \frac{5}{8}$ of the integral of the coarse bubble.
Hence, when using a hierarchical basis, since the piecewise linear basis
functions are prolonged exactly we can obtain a prolongation that satisfies
\eqref{eqn:decomposition-exterior} by simply multiplying the coefficients
of the fine grid bubble functions by $8/5$. After this scaling, the local
correction $w_T$ is computed as described above.
For a nodal basis, a change of basis to the hierarchical basis should be
performed.

This modification of the prolongation operator is crucial for the solver to work with the
$[\PoneFB]^3\mathrm{-}\Pzero$ element. We demonstrate this by showing the
residual of the outer flexible GMRES iteration for the linear solve in the
first Newton step at $\Re=10$ for a lid-driven cavity problem (see section
\ref{ssec:threedimensionalexperiments} for details) in Table~\ref{tab:ksp-failure-without-scaling}.
Without modifying the prolongation of the facet bubbles, we observe no convergence.
\begin{table}[htbp]
    \centering
    \begin{tabular}{r|c|c}
        Iteration & Residual with bubble scaling & Residual without bubble scaling \\\hline
        0         &         $ 3.499\hphantom{{}\times 10^{+0}}$ & $ 3.499$            \\
        1         &         $ 1.554\times 10^{-2}$ & $ 3.499 $            \\
        2         &         $ 1.716\times 10^{-3}$ & $ 3.499 $            \\
        3         &         $ 1.821\times 10^{-4}$ & $ 3.496 $            \\
        4         &         $ 1.651\times 10^{-5}$ & $ 3.495 $            \\
    \end{tabular} 
    \caption{Residual of the outer flexible GMRES solver when employing the
    $[\PoneFB]^3\mathrm{-}\Pzero$ element. It is necessary to modify the prolongation operator
    to achieve convergence with this element.}
    \label{tab:ksp-failure-without-scaling}
\end{table}

Lastly, we consider the $[\PtwoFB]^3-\Pzero$ element. While it is also
non-nested, it turns out that the interpolation is exact \emph{on the facets}
of each coarse cell and hence flux preserving. To see this, observe
that the cubic facet bubble function is only quadratic on the newly introduced
edges of a regularly refined facet, as they are parallel to the edges of the
coarse facet and therefore one of the barycentric coordinates is constant. The coarse bubble function
is therefore prolonged exactly. This means that the $[\PtwoFB]^3-\Pzero$
element can be used with the Sch\"oberl prolongation operator \eqref{eqn:schoberlprolong},
without the modifications necessary for $[\PoneFB]^3\mathrm{-}\Pzero$ described
above. However, in our preliminary numerical experiments the simpler
prolongation was outweighed by the cost of the larger number of degrees
of degrees of freedom, and hence we use $[\PoneFB]^3\mathrm{-}\Pzero$ for
the numerical experiments in section \ref{sec:numerical}.

\begin{remark}
Only the prolongation is modified; as in Benzi \& Olshanskii \cite{benzi2006}, the natural operations are used for
restriction and injection.
\end{remark}

\subsection{The advection terms} \label{sec:advection}
So far we have neglected the terms arising from the linearization of the
advection term. Applying a Newton linearization, \eqref{eqn:bilinearstokes}
becomes: find $u \in V_{h,0}$ such that
\begin{equation} \label{eqn:advectnewton}
(2\nu\eps{u}, \nabla v) + (w \cdot \nabla u, v) + (u \cdot \nabla w, v) + \gamma (\Pq \nabla \cdot u, \nabla \cdot
v) = (f, v)
\end{equation}
for all $v \in V_{h,0}$, while the Picard linearization yields:
find $u \in V_{h,0}$ such that
\begin{equation} \label{eqn:advectpicard}
(2\nu\eps{u}, \nabla v) + (w \cdot \nabla u, v) + \gamma (\Pq \nabla \cdot u, \nabla \cdot
v) = (f, v)
\end{equation}
for all $v \in V_{h,0}$. The Picard linearization is easier to solve but
sacrifices quadratic convergence of the nonlinear solver. Several authors have
reported success with geometric multigrid for scalar analogues of
\eqref{eqn:advectpicard} without the grad-div term, using a combination of
line/plane relaxation and SUPG stabilization
\cite{ramage1999,olshanskii2004,wu2006}. Olshanskii and Benzi
\cite{olshanskii2008} and Elman et al.~\cite{elman2003} apply preconditioners built on the Picard
linearization \eqref{eqn:advectpicard} to the Newton linearization
\eqref{eqn:advectnewton}, with good results.

Numerical experiments indicated that the additive star iteration alone
was not effective as a relaxation method for
\eqref{eqn:advectnewton}. (Benzi and Olshanskii \cite{benzi2006} used
a multiplicative star iteration with multiple directional sweeps, but
we wished to avoid this as its performance varies with the core count
in parallel.) We investigated the multiplicative composition of
additive star iterations and plane smoothers, and while this led to a
successful multigrid cycle, the plane smoothers were quite expensive
(involving many 2D solves) and were also difficult to parallelize on
arbitrary unstructured grids where the parallel decomposition does not
divide into planes. While the additive star iteration alone is not
effective as a relaxation for \eqref{eqn:advectnewton}, we found that a
few iterations of GMRES preconditioned by the additive star iteration
was surprisingly effective as a relaxation method, even for low viscosities. This point merits
further analysis and will be considered in future work. This
relaxation method also has the advantage that it is easy to
parallelize, with convergence independent of the parallel
decomposition.

\section{Numerical results} \label{sec:numerical}

\subsection{Algorithm details}
\begin{figure}[htbp]
  \centering
  \resizebox{\textwidth}{!}{\begin{tikzpicture}[%
  every node/.style={draw=black, thick, anchor=west},
  grow via three points={one child at (0.0,-0.7) and
  two children at (0.0,-0.7) and (0.0,-1.4)},
  edge from parent path={(\tikzparentnode.210) |- (\tikzchildnode.west)}]
  \node {Continuation}
    child { node {Newton solver with line search}
      child { node {Krylov solver (FGMRES)}
        child { node {Block preconditioner}
          child { node {Approximate Schur complement inverse}}
          child { node {F-cycle on augmented momentum block}
              child { node {Coarse grid solver}
                child { node {LU factorization}}
              }
              child [missing] {}
              child { node {Prolongation operator}
                child { node {Local solves over coarse cells}}
              }
              child [missing] {}
              child { node {Relaxation}
                child { node {GMRES}
                  child { node {Additive star iteration}}
                }
              }
          }
        }
      }
    };
\end{tikzpicture}}
  \caption{An outline of the algorithm for solving \eqref{eqn:ns}.}
  \label{fig:solver}
\end{figure}
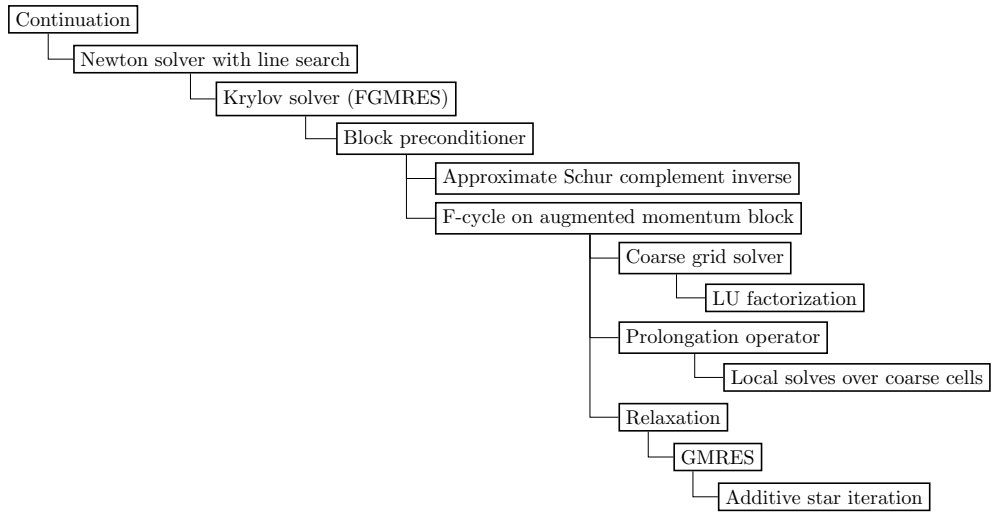

A graphical representation of the entire algorithm is shown in Figure
\ref{fig:solver}.
We employ simple continuation in Reynolds number as a globalization
device, as Newton's method is not globally convergent. Newton's method is globalized
with the $L^2$ line search algorithm of PETSc \cite{brune2015}.

We use flexible
GMRES \cite{saad1993} as the outermost solver for the linearized Newton system,
as we employ GMRES in the multigrid relaxation.
If the pressure is only defined up to a
constant, then the appropriate nullspace is passed to the Krylov solver and the
solution is orthogonalized against the nullspace at every iteration.
We use the full block factorization
preconditioner
\begin{equation}
P^{-1} =
\left( \begin{array}{cc}
I   & -\tilde{A}_\gamma^{-1} B^T \\
0 & I \\
\end{array} \right)
\left( \begin{array}{cc}
\tilde{A}_\gamma^{-1}  & 0 \\
0 & \tilde{S}^{-1} \\
\end{array} \right)
\left( \begin{array}{cc}
I   & 0 \\
-B\tilde{A}_\gamma^{-1} & I \\
\end{array} \right)
\end{equation}
with approximate inner solves $\tilde{A}_\gamma^{-1}$ and $\tilde{S}^{-1}$ for the
augmented momentum block and the Schur complement respectively. The diagonal,
upper and lower triangular variants described in \cite{murphy2000,ipsen2001}
also converge well, but these took longer runtimes in preliminary experiments.

We use one application of a full multigrid cycle \cite[Figure
1.2]{brandt2011} using the components described in section
\ref{sec:solvingmomentum} for $\tilde{A}_\gamma^{-1}$. The problem on each level
is constructed by rediscretization; fine grid functions, such as the current
iterate in the Newton scheme, are transferred to the coarse levels via
injection. On each level the SUPG stabilization is performed with parameters
corresponding to the mesh in question. For each relaxation
sweep we perform 6 (in 2D) or 10 (in 3D) GMRES iterations preconditioned by the
additive star iteration; at lower Reynolds numbers this can be reduced, but we
found that these expensive smoothers represented the optimal tradeoff between
inner and outer work at higher Reynolds numbers. The problem on the coarsest level is solved with the SuperLU\_DIST sparse
direct solver \cite{li1999,li2003b}. For scalability, the coarse grid solve is
agglomerated onto a single compute node using PETSc's telescoping facility
\cite{may2016}. As all inner solvers are additive, the convergence of the solver
is independent of the parallel decomposition (up to roundoff).

\subsection{Software implementation}
The solver proposed in the previous section is complex, and relies heavily on
PETSc's capability for the arbitrarily nested composition of solvers
\cite{brown2012}. For the implementation of local patch solves, we have
developed a new subspace correction preconditioner for PETSc that relies on the
\texttt{DMPlex} unstructured mesh component \cite{knepley2005,knepley2009} for
topological subspace construction and provides an extensible callback interface
that allows for the very general specification of additive Schwarz methods. A
detailed description of this preconditioner is in preparation.

\subsection{Solver verification with the method of manufactured solutions}
In order to verify the implementation and the convergence of the
$[\PoneFB]^3\mathrm{-}\Pzero$ element we employ the method of manufactured
solution.
We start by considering the pressure and velocity field proposed in \cite{shih1989}, which is rescaled 
to the $[0, 2]^2$ square.
This results in $u=(u_1, u_2)$ with
\begin{equation}
    \begin{aligned}
        u_1(x, y) =&\  \frac{1}{4} (x-2)^2 x^2 y \left(y^2-2\right) \\
        u_2(x, y) =&\  -\frac{1}{4} x \left(x^2-3 x+2\right) y^2 \left(y^2-4\right) \\
        \tilde p(x, y) =&\  \frac{x y \left(3 x^4-15 x^3+10 x^2 y^2-30 x \left(y^2-2\right)+20 \left(y^2-2\right)\right)}{5 \Re}\\
                    &\  -\frac{1}{128} (x-2)^4 x^4 y^2 \left(y^4-2 y^2+8\right)\\
        p(x, y)= &\ \tilde p -\frac14\int_{[0, 2]^2} \tilde p(x,y)\dx = \tilde p  +\frac{1408}{33075} - \frac{8}{5\Re}.
    \end{aligned}
\end{equation}
As we are primarily interested in the three dimensional case, we extend the
vector field into the $z$ dimension via $u(x, y, z) = (u_1(x, y), u_2(x, y), 0)$.
The pressure remains the same as in two dimensions.

\begin{figure}[htbp]
  \centering
  \pgfplotstableread[col sep=comma, row sep=\\]{%
    h,error_v,error_p\\
    0.125,0.0469998769846097,4.879481060699821\\
    0.0625,0.013324555990202134,1.7148247471726807\\
    0.03125,0.0034906304895967115,0.6896606226331777\\
    0.015625,0.0008795398151003766,0.31243679597680674\\}\reonegammatenthousand

  \pgfplotstableread[col sep=comma, row sep=\\]{%
    h,error_v,error_p\\
    0.125,0.04699987762213687,4.8794810191317\\
    0.0625,0.013324556679585035,1.7148196381504333\\
    0.03125,0.0034906116893429814,0.6896557435605192\\
    0.015625,0.000879529966251674,0.3124295483848757\\}\reonegammaone
  \pgfplotstableread[col sep=comma, row sep=\\]{%
    h,error_v,error_p\\
    0.125,0.10417962511049522,0.06313533525600275\\
    0.0625,0.03952463294902316,0.026975093770353203\\
    0.03125,0.012085756491489108,0.012039184365974821\\
    0.015625,0.003257861159774448,0.0057070900679878364\\}\retwohundredgammatenthousand

  \pgfplotstableread[col sep=comma, row sep=\\]{%
    h,error_v,error_p\\
    0.125,0.1041796371676739, 0.0631353304155226\\
    0.0625,0.03952461999640074,0.026975118503031525\\
    0.03125,0.012085808316116402,0.01203917628623866\\
    0.015625,0.003257814558224604,0.005707086017297249\\}\retwohundredgammaone
  \pgfplotstableread[col sep=comma, row sep=\\]{%
    h,error_v,error_p\\
    0.125,0.2502014617647885,0.0822492048847913\\
    0.0625,0.14684788158100423,0.043708556176756605\\
    0.03125,0.057281809559438804,0.01797688614863548\\
    0.015625,0.017045784991612873,0.0068460314194125835\\}\refivehundredgammatenthousand

  \pgfplotstableread[col sep=comma, row sep=\\]{%
    h,error_v,error_p\\
    0.125,0.25020149507893197,0.08224921243261697\\
    0.0625,0.14684789993264036,0.04370852728148409\\
    0.03125,0.05728182350093682,0.017976847105105765\\
    0.015625,0.01704621605013802,0.006846097701282368\\}\refivehundredgammaone
  \begin{tikzpicture}
    \begin{loglogaxis}[
      width=0.5\linewidth,
      xlabel={$h$},
      ylabel near ticks,
      ylabel={$L^2$-error},
      log basis x=2,
      title={$\Re = 1$},
      ]
      \addplot+[thick, mark=o, mark options={fill=none}, mark size=4pt] table[x=h, y=error_v] {\reonegammatenthousand};
      \addplot+[thick, mark=square, mark options={fill=none}, mark size=2pt] table[x=h, y=error_v] {\reonegammaone};
      \addplot+[black, dotted, mark=none] table[x=h, y expr={0.5*\thisrow{h}^2}]
      {\reonegammaone};

      \addplot+[thick, mark=pentagon, mark options={fill=none}, mark size=4pt] table[x=h, y=error_p] {\reonegammatenthousand};
      \addplot+[thick, mark=triangle, mark options={fill=none}, mark size=2pt] table[x=h, y=error_p] {\reonegammaone};
      \addplot+[black, dashed, mark=none] table[x=h, y expr={50*\thisrow{h}}]
      {\reonegammaone};
    \end{loglogaxis}
  \end{tikzpicture}
  \begin{tikzpicture}
    \begin{loglogaxis}[
      width=0.5\linewidth,
      xlabel={$h$},
      ylabel near ticks,
      ylabel={$L^2$-error},
      log basis x=2,
      title={$\Re = 200$}
      ]
      \addplot+[thick, mark=o, mark options={fill=none}, mark size=4pt] table[x=h, y=error_v] {\retwohundredgammatenthousand};
      \addplot+[thick, mark=square, mark options={fill=none}, mark size=2pt] table[x=h, y=error_v] {\retwohundredgammaone};
      \addplot+[black, dotted, mark=none] table[x=h, y expr={20*\thisrow{h}^2}]
      {\retwohundredgammaone};

      \addplot+[thick, mark=pentagon, mark options={fill=none}, mark size=4pt] table[x=h, y=error_p] {\retwohundredgammatenthousand};
      \addplot+[thick, mark=triangle, mark options={fill=none}, mark size=2pt] table[x=h, y=error_p] {\retwohundredgammaone};
      \addplot+[black, dashed, mark=none] table[x=h, y expr={0.3*\thisrow{h}}]
      {\reonegammaone};
    \end{loglogaxis}
  \end{tikzpicture}
  \begin{tikzpicture}
    \begin{loglogaxis}[
      width=0.5\linewidth,
      xlabel={$h$},
      ylabel near ticks,
      ylabel={$L^2$-error},
      log basis x=2,
      title={$\Re = 500$}, 
      legend cell align=left,
      legend style={at={(1.425, 0.5)}, anchor=west},
      legend columns=1,
      legend style={/tikz/every even column/.append style={column
          sep=0.5cm}}     
      ]
      \addplot+[thick, mark=o, mark options={fill=none}, mark size=4pt] table[x=h, y=error_v] {\refivehundredgammatenthousand};
      \addlegendentry{$\|v - v_h\|_{L^2}$, $\gamma = 10^4$}
      \addplot+[thick, mark=square, mark options={fill=none}, mark size=2pt] table[x=h, y=error_v] {\refivehundredgammaone};
      \addlegendentry{$\|v - v_h\|_{L^2}$, $\gamma = 1$}
      \addplot+[black, dotted, mark=none] table[x=h, y expr={100*\thisrow{h}^2}]
      {\reonegammaone};
      \addlegendentry{$h^2$};

      \addplot+[thick, mark=pentagon, mark options={fill=none}, mark size=4pt] table[x=h, y=error_p] {\refivehundredgammatenthousand};
      \addlegendentry{$\|p - p_h\|_{L^2}$, $\gamma = 10^4$}
      \addplot+[thick, mark=triangle, mark options={fill=none}, mark size=2pt] table[x=h, y=error_p] {\refivehundredgammaone};
      \addlegendentry{$\|p - p_h\|_{L^2}$, $\gamma = 1$}
      \addplot+[black, dashed, mark=none] table[x=h, y expr={0.3*\thisrow{h}}]
      {\reonegammaone};
      \addlegendentry{$h$};
    \end{loglogaxis}
    \path[use as bounding box] (current bounding box.south west)
    rectangle ($(current bounding box.north east) + (0.4605cm, 0)$);
  \end{tikzpicture}
  \caption{Convergence of the computed velocity and pressure field as the
    mesh is refined for a 3D lid-driven cavity test problem.}
  \label{fig:mms}
\end{figure}
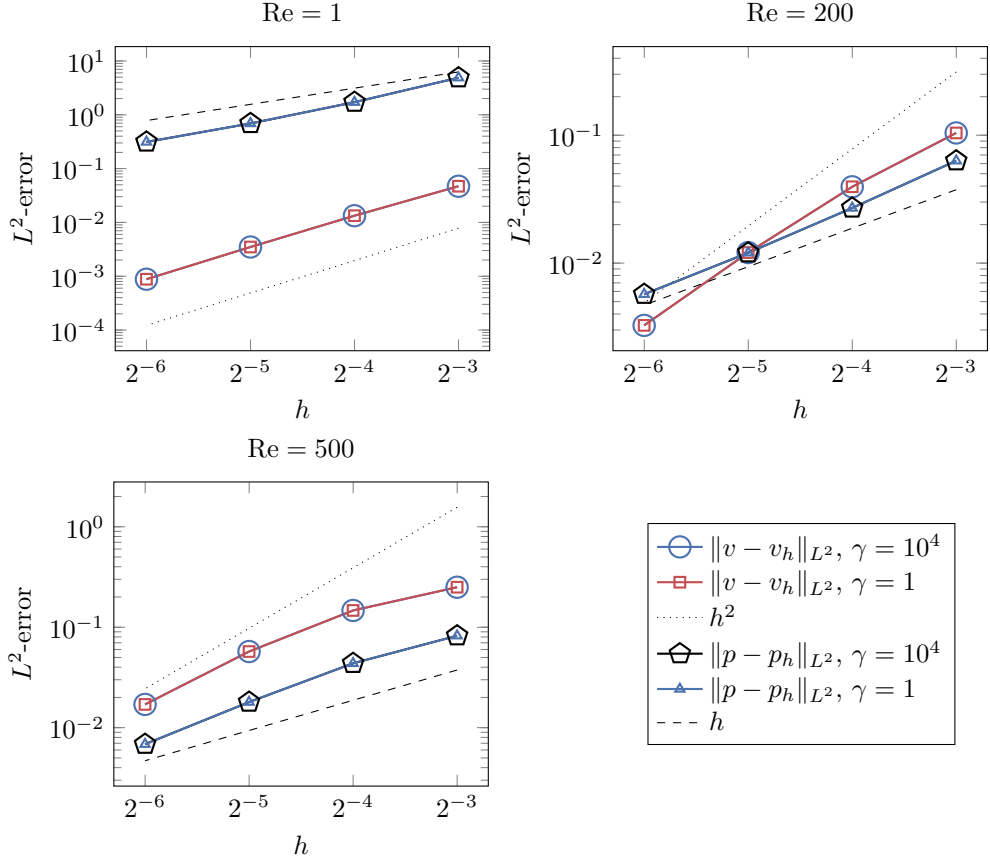

To demonstrate that the error convergence is independent of $\gamma$, we run
the solver for values $\gamma=1$ and $\gamma=10^4$. Figure~\ref{fig:mms} shows
the error between the computed velocity and pressure and their known
analytical solutions for $\Re=1$, $\Re=200$ and $\Re=500$. Due to the DG0
discretization we expect, and see, first order convergence of the pressure.
Without stabilization, we expect second order convergence for the
velocity field; however, due to the presence of the SUPG stabilization
this is reduced to linear convergence for coarse meshes. Once the mesh
is fine enough so that $h^{-1} \gtrsim \Re$, second order convergence is recovered.

\subsection{Two-dimensional experiments}
We consider two representative benchmark problems: the regularized lid-driven
cavity and backward-facing step problems, fully described in 
\cite[examples 8.1.2 and 8.1.3]{elman2014}. For each experiment, we fix a coarse
grid and vary the number of refinements to vary the size of the problem under
consideration; all refinements are used in the multigrid iteration, to ensure
that the convergence does not deteriorate as more levels are employed. We
employ the $[\Ptwo]^2\mathrm{-}\Pzero$ element for all two dimensional experiments. 
To investigate the performance of the solver with Reynolds number, the problem is
first solved for $\Re = 10$, then $\Re = 100$, and then in steps of $100$ until $\Re =
10000$, with the solution for the previous value of \Re used as initial guess
for the next. The Stokes equations are solved using a standard geometric
multigrid algorithm with the pressure mass matrix as Schur complement
approximation and point-block SOR as a smoother to provide the initial guess
used at $\Re = 10$. The augmented Lagrangian parameter is set to $\gamma =
10^4$ in these and all subsequent experiments.

The linear solves are terminated with an absolute tolerance of $10^{-10}$ in
the $\ell_2$-norm and a relative tolerance of $10^{-6}$. The nonlinear solves
are terminated with an absolute tolerance of $10^{-8}$ and a relative tolerance
of $10^{-10}$. As each outer iteration of the Krylov method does a
\emph{fixed} amount of work (i.e.~all subproblems are solved with a fixed
number of iterations, not to a specified tolerance), the solver scales well
with mesh size and Reynolds number if the iteration counts remain approximately
constant.

For comparison, we solve the same problems using the reference implementations
of the PCD and LSC preconditioners in version 3.5 of IFISS \cite{elman2014b}, up
to $\Re = 1000$, as IFISS does not employ stabilization of the advection term.
For both of these preconditioners we use the variant that takes corrections for
the boundary conditions into account and we solve the inner problems in the
Schur complement approximation using an algebraic multigrid solver. We employ
the hybrid strategy suggested by \cite[p. 391]{elman2014} that uses a single
sweep of ILU(0) on the finest level and two iterations of point-damped Jacobi
for pre- and post-smoothing on all coarsened levels. A relative tolerance of
$10^{-6}$ is set for the Krylov solver and an absolute tolerance of $10^{-8}$
for the Newton solver.

We begin by considering the regularized lid-driven cavity problem. The coarse
grid used is the $16 \times 16$ grid of triangles of negative slope. The
results are shown in Table~\ref{tab:ourldc}; we observe only very mild
iteration growth from $\Re=10$ to $\Re=10000$ with the performance improving as
more refinements are taken. Iteration counts using the PCD and LSC
preconditioners are shown in Table~\ref{tab:pcdlscldc}. For both PCD and LSC
iteration counts increase substantially from $\Re = 10$ to $\Re = 1000$.
\begin{table}[htbp]
\centering
\begin{tabular}{cc|ccccc}
\toprule
\# refinements & \# degrees of freedom & \multicolumn{5}{c}{Reynolds number} \\
  && 10 & 100 & 1000 & 5000 & 10000 \\
\midrule
1 & $1.0 \times 10^4$ & 2.50 &  4.33 &  6.00 & 8.00 & 14.00 \\
2 & $4.1 \times 10^4$ & 2.50 &  3.33 &  6.67 & 8.50 & 10.00 \\
3 & $1.6 \times 10^5$ & 2.50 &  3.00 &  5.67 & 8.50 & 9.00 \\
4 & $6.6 \times 10^5$ & 2.50 &  2.67 &  5.00 & 8.00 & 8.50 \\
\bottomrule
\end{tabular}
\caption{Average number of outer Krylov iterations per Newton step for the
2D regularized lid-driven cavity problem.}
\label{tab:ourldc}
\end{table}

\begin{table}[htbp]
\centering
\begin{tabular}{cc|c@{\hspace{9pt}}c@{\hspace{9pt}}c}
\toprule
$1/h$ & \# degrees of freedom & \multicolumn{3}{c}{Reynolds number} \\
  && 10 & 100 & 1000 \\
\midrule
$2^4$ & $8.34 \times 10^2$ & 22.0/21.5 & 40.4/48.7 & 103.3/130.7 \\
$2^5$ & $3.20 \times 10^3$ & 23.0/22.0 & 41.3/52.7 & 137.7/185.3 \\
$2^6$ & $1.25 \times 10^4$ & 24.5/22.5 & 42.0/49.3 & 157.0/205.7 \\
$2^7$ & $4.97 \times 10^4$ & 25.5/21.0 & 42.7/43.3 & 149.0/207.3 \\
$2^8$ & $1.98 \times 10^5$ & 26.0/23.0 & 44.0/38.0 & 137.0/180.0 \\
\bottomrule
\end{tabular}
\caption{Average number of outer Krylov iterations per Newton step for the
2D regularized lid-driven cavity problem with PCD/LSC preconditioner.}
\label{tab:pcdlscldc}
\end{table}

For the backward-facing step we observe that the performance is dependent on
the resolution of the coarse grid. We consider two experiments, one starting
with a coarse grid consisting of 6941 vertices and 13880 elements (labeled A)
and one consisting of 30322 vertices and 60642 elements (labeled B). Both
unstructured triangular meshes were generated with Gmsh \cite{geuzaine2009}.
For mesh A, we observe that the iteration counts for large Reynolds
numbers show the solver degrades somewhat as the mesh is
refined, see Table~\ref{tab:ourbfs}. Using the finer coarse grid B alleviates
this problem. The bottom half of Table~\ref{tab:ourbfs} shows that iteration counts only
approximately double as we increase from $\Re=10$ to $\Re=10000$.
\begin{table}[htbp] 
\centering
\begin{tabular}{cc|ccccc}
\toprule
\# refinements & \# degrees of freedom & \multicolumn{5}{c}{Reynolds number} \\
               && 10 & 100 & 1000 &  5000 & 10000\\
\midrule
\multicolumn{7}{c}{coarse grid A}\\
\midrule
1 & $2.7 \times 10^5$ & 3.00 & 4.00 & 5.50 & 12.00 & 26.50\\
2 & $1.1 \times 10^6$ & 2.67 & 4.50 & 5.50 & 11.50 & 31.50\\
3 & $4.3 \times 10^6$ & 4.00 & 7.00 & 6.00 & 14.00 & 21.00\\
\midrule
\multicolumn{7}{c}{coarse grid B}\\
\midrule
1 & $1.2 \times 10^6$ & 2.67 & 3.75 & 5.00 & 7.00 & 9.00\\
2 & $4.8 \times 10^6$ & 4.00 & 3.75 & 5.00 & 6.50 & 7.00\\
3 & $1.9 \times 10^7$ & 3.67 & 5.75 & 5.50 & 5.00 & 5.50\\
\bottomrule
\end{tabular}
\caption{Average number of outer Krylov iterations per Newton step for the
2D backward-facing step problem for two different coarse grids.}
\label{tab:ourbfs}
\end{table}

The results for PCD and LSC on the backwards-facing step are shown in Table
\ref{tab:pcdlscbfs}. The iteration counts approximately treble as we increase
from $\Re = 10$ to $\Re = 1000$.
\begin{table}[htbp]
\centering
\begin{tabular}{cc|c@{\hspace{10pt}}c@{\hspace{10pt}}c}
\toprule
$1/h$ & \# degrees of freedom & \multicolumn{3}{c}{Reynolds number} \\
      && 10 & 100 & 1000 \\
\midrule
$2^4$ & $3.94 \times 10^3$ & 23.0/29.0 & 32.5/47.5 & \texttt{NaNF}/\texttt{NaNF} \\
$2^5$ & $1.52 \times 10^4$ & 23.5/26.0 & 31.0/45.0 & 221.3/329.0                 \\
$2^6$ & $5.96 \times 10^4$ & 23.5/25.5 & 30.5/42.8 & 122.7/225.7                 \\
$2^7$ & $2.36 \times 10^5$ & 23.5/25.5 & 30.0/40.8 & 85.3/161.3                  \\
$2^8$ & $9.38 \times 10^5$ & 23.5/27.0 & 30.0/40.0 & 78.3/128.0                  \\
\bottomrule
\end{tabular}
\caption{Average number of outer Krylov iterations per Newton step for the 2D
backwards-facing step problem with PCD/LSC preconditioner. \texttt{NaNF}
denotes failure due to NaNs occurring in the solve for the velocity block.}
\label{tab:pcdlscbfs}
\end{table}

\subsubsection{Runtime comparison to SIMPLE}
\label{sec:comparison-simple}

To give some measure of the runtime of the solver, we compare it to an
implementation of SIMPLE~\citep[Section 6.7]{patankar1980} in the same software framework. We select the
lid-driven cavity in two dimensions with three refinements ($1.6 \times 10^5$ degrees of freedom) as a representative
problem. The SIMPLE preconditioner is given by
\begin{equation}
P_\text{SIMPLE}^{-1} =
\left( \begin{array}{cc}
I   & -\operatorname{diag}(A)^{-1} B^T \\
0 & I \\
\end{array} \right)
\left( \begin{array}{cc}
\tilde{A}^{-1}  & 0 \\
0 & \tilde{S}_\text{SIMPLE}^{-1} \\
\end{array} \right)
\left( \begin{array}{cc}
I   & 0 \\
-B\tilde{A}^{-1} & I \\
\end{array} \right),
\end{equation}
where
\begin{equation}
  \tilde{S}_\text{SIMPLE} = -B^T \operatorname{diag}(A)^{-1} B
\end{equation}
and no grad-div augmentation is employed. $\tilde{A}^{-1}$ is approximated
by one full multigrid cycle of the ML
algebraic multigrid library~\citep{gee2006}; $\tilde{S}_\text{SIMPLE}^{-1}$ is
approximated with one V cycle of ML\footnote{For fairness, we do not use exact inner
solves, since our solver also does not use exact inner solves. Of the
algebraic multigrid libraries available in PETSc, ML performed the best with default settings.}.

\begin{table}[htbp]
\centering
\begin{tabular}{c|c@{\hspace{10pt}}c|c@{\hspace{10pt}}c}
\toprule
Reynolds number & \multicolumn{2}{c|}{Augmented Lagrangian} & \multicolumn{2}{c}{SIMPLE}                 \\
                & Total iterations                          & Time (min) & Total iterations & Time (min) \\
  \midrule
  10            & 4                                         & 0.21       & 515              & 1.06       \\ 
  50            & 6                                         & 0.30       & 741              & 1.57       \\ 
  100           & 8                                         & 0.38       & 979              & 2.01       \\ 
  150           & 10                                        & 0.46       & 1111             & 2.27       \\ 
  200           & 10                                        & 0.46       & 1185             & 2.48       \\ 
\bottomrule
\end{tabular}
\caption{Iteration count and runtime comparison against the SIMPLE preconditioner.}
\label{tab:simplealmgcompare}
\end{table}

The results for several continuation steps are shown in Table
\ref{tab:simplealmgcompare}. The computations were performed in serial. Each SIMPLE iteration is approximately 22--26 times
faster than an augmented Lagrangian iteration, but the lower cost per iteration
is outweighed by the greater number of iterations required.

\subsection{Three-dimensional experiments}
\label{ssec:threedimensionalexperiments}
The lid-driven cavity and backward-facing step problems can both
be extended to three dimensions in a natural way.
For the lid-driven cavity, we consider the cube $\Omega =[0,2]^3$ with no-slip boundary conditions on all sides apart from the top boundary $\{y=2\}$.
On the top boundary we enforce $u(x, y, z)=(x^2(2-x)^2z^2(2-z)^2, 0, 0)^T$.
The three dimensional backwards-facing step is given by $\Omega=\left(\left([0, 10]\times [1, 2]\right)\cup \left([1, 10] \times [0, 1]\right)\right)  \times [0, 1]$.
We enforce the inflow condition $u(x, y, z) = (4(2-y)(y-1)z(1-z), 0, 0)^T$ on the
left boundary $\{x=0\}$, a natural outflow condition on the right boundary
$\{x=10\}$ and no-slip boundary conditions on the
remaining boundaries.

\begin{figure}[htbp]
  \centering
    \includegraphics[width=0.6\textwidth]{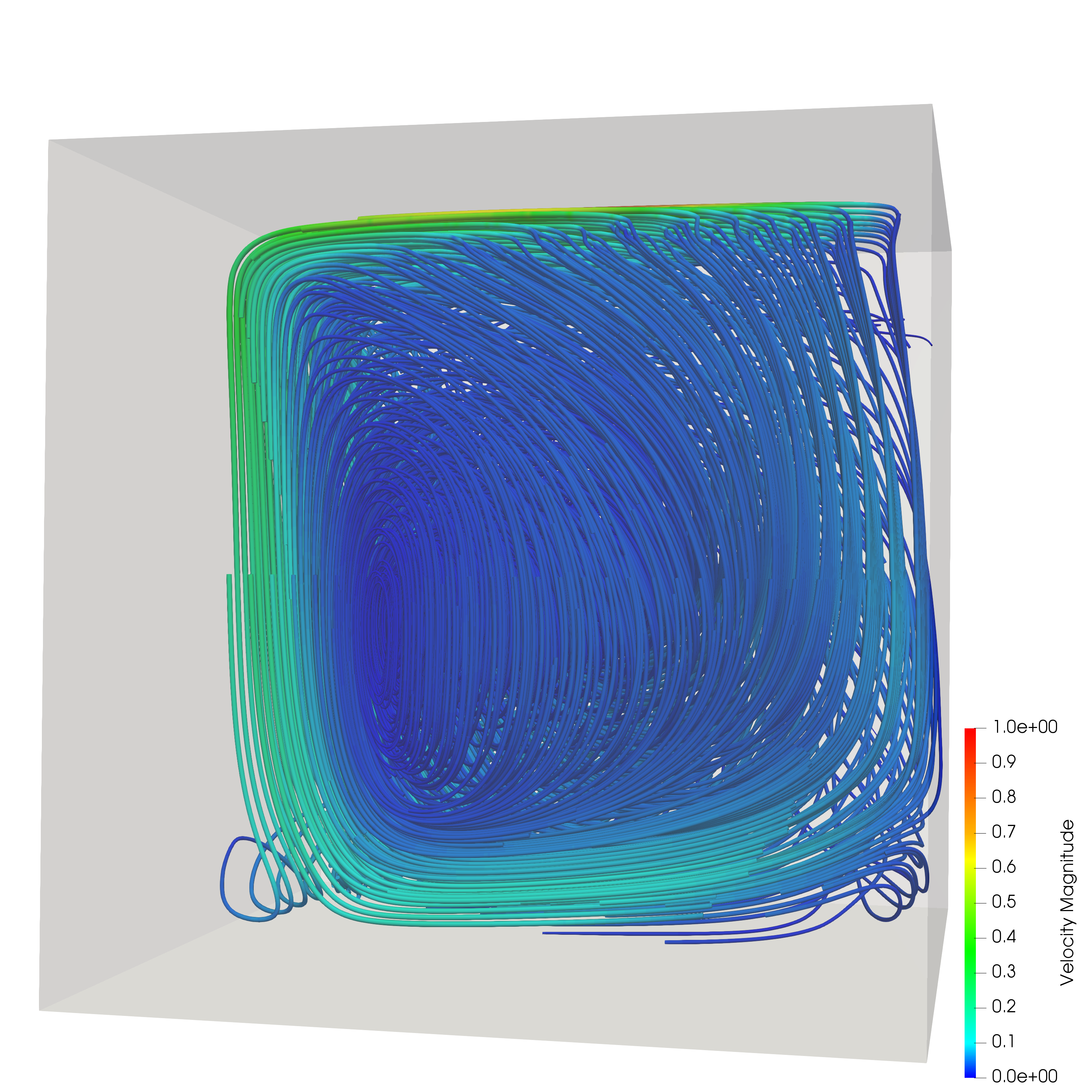}
  \caption{Interior view of the streamtubes of the 3D lid-driven cavity at $\Re = 5000$. The color
  denotes speed.}
  \label{fig:bfs}
\end{figure}

\begin{figure}[htbp]
  \centering
    \includegraphics[width=\textwidth]{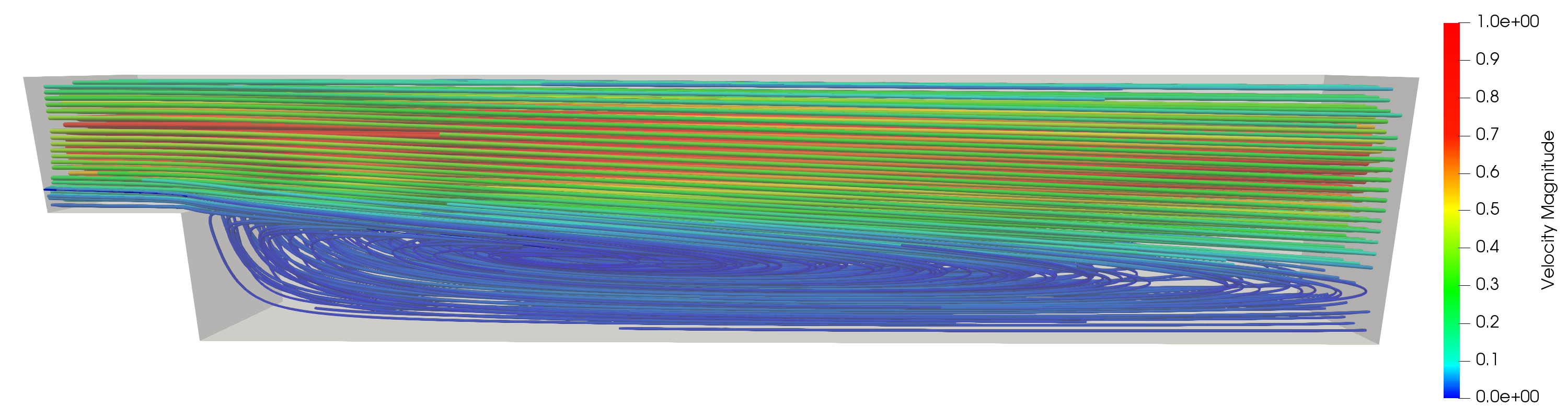}
  \caption{Interior view of the streamtubes of the 3D backwards-facing step at $\Re = 1000$. The color
  denotes speed.}
  \label{fig:ldc}
\end{figure}

Two aspects of the solver were modified compared to the version used in two
dimensions. Firstly, we observe that reducing the size of the SUPG
stabilization by a factor of $1/20$ improves convergence significantly.
Secondly, the relative tolerance for the linear solver was relaxed to
$10^{-5}$ and the absolute tolerance for the linear and the nonlinear solver
was relaxed to $10^{-8}$, to save computational time. The three-dimensional
experiments were both run for $[\PoneFB]^3\mathrm{-}\Pzero$ discretizations of
up to one billion degrees of freedom on ARCHER, the UK national supercomputer.
The runs were terminated at $\Re = 5000$ due to budgetary constraints.
Images of the solutions are shown in Figures \ref{fig:bfs} and \ref{fig:ldc}.
\begin{table}[htbp]
\centering
\begin{tabular}{cc|ccccc}
\toprule
\# refinements & \# degrees of freedom & \multicolumn{5}{c}{Reynolds number} \\
  && 10 & 100 & 1000 & 2500 & 5000 \\
\midrule
1 & $2.1 \times 10^6$ & 4.50 & 4.00 & 5.00 & 4.50 & 4.00 \\
2 & $1.7 \times 10^7$ & 4.50 & 4.33 & 4.50 & 4.00 & 4.00 \\
3 & $1.3 \times 10^8$ & 4.50 & 4.33 & 4.00 & 3.50 & 7.00 \\
4 & $1.1 \times 10^9$ & 4.50 & 3.66 & 3.00 & 5.00 & 5.00 \\
\bottomrule
\end{tabular}
\caption{Average number of outer Krylov iterations per Newton step for the
3D regularized lid-driven cavity problem.}
\label{tab:ourldc3d}
\end{table}

\begin{table}[htbp]
\centering
\begin{tabular}{cc|ccccc}
\toprule
\# refinements & \# degrees of freedom & \multicolumn{5}{c}{Reynolds number} \\
               && 10 & 100 & 1000 & 2500& 5000 \\
\midrule
1 & $2.1 \times 10^6$ & 4.50 & 4.00 & 4.00 & 4.50 & 7.50  \\
2 & $1.7 \times 10^7$ & 5.00 & 4.00 & 3.33 & 4.00 & 10.00 \\
3 & $1.3 \times 10^8$ & 6.50 & 4.50 & 3.50 & 3.00 & 8.00  \\
4 & $1.0 \times 10^9$ & 7.50 & 3.50 & 2.50 & 3.00 & 6.00  \\
\bottomrule
\end{tabular}
\caption{Average number of outer Krylov iterations per Newton step for the
3D backwards-facing step problem.}
\label{tab:ourbfs3d}
\end{table}
As for the two-dimensional case, we see only very little variation of the
iteration counts with Reynolds number over this range.

To stress the solver further, the lid-driven cavity with 2 refinements
($1.7 \times 10^7$ degrees of freedom)
was run until failure. Iteration counts remain stable until $\Re = 7000$,
then begin to increase, with eventual failure of convergence at $\Re = 7700$.


\subsubsection{Computational performance}
\label{sec:comp-perf}

\begin{figure}[htbp]
  \centering
  \pgfplotstableread[col sep=comma, row sep=\\]{%
    Cores,Time,Dofs\\
    48,5.3e1,2134839\\
    384,6.9e1,16936779\\
    3072,6.1e1,134930451\\
    24576,6.6e1,1077196323\\
  }\ldctable
  \pgfplotstableread[col sep=comma, row sep=\\]{%
    Cores,Time,Dofs\\
    48,5.5e1,2107839\\
    384,6.9e1,16534263\\
    3072,5.7e1,130973115\\
    24576,6.9e1,1042606515\\
  }\bfstable
  \begin{subfigure}[l]{0.475\textwidth}
    \begin{tikzpicture}
      \begin{semilogxaxis}[
        width=0.925\linewidth,
        height=0.8\linewidth,
        log basis x=2,
        ymin=0,
        ymax=80,
        xtick=data,
        xticklabels from table={\ldctable}{Cores},
        extra x ticks={48, 384, 3072, 24576},
        extra x tick labels={$[2.13]$, $[16.9]$,$[135]$,$[1077]$},
        extra x tick style={tick label style={yshift=-2.5ex}},
        xlabel={Cores\\{}[DoFs $\times 10^6$]},
        xlabel style={align=center, style={yshift=-2ex}},
        ylabel near ticks,
        ylabel style={align=center, text width=4cm},
        ylabel={Time to solution over all continuation steps [min]},
        ]
        \addplot+ table[x=Cores,y=Time] {\ldctable};
      \end{semilogxaxis}
    \end{tikzpicture}
    \caption{Weak scaling of the three-dimensional lid-driven cavity.}
    \label{fig:ldc3d-scaling}
  \end{subfigure}
  \hspace{1em}
  \begin{subfigure}[r]{0.475\textwidth}
    \begin{tikzpicture}
      \begin{semilogxaxis}[
        width=0.925\linewidth,
        height=0.8\linewidth,
        log basis x=2,
        ymin=0,
        ymax=80,
        xtick=data,
        xticklabels from table={\bfstable}{Cores},
        extra x ticks={48, 384, 3072, 24576},
        extra x tick labels={$[2.11]$, $[16.5]$,$[131]$,$[1043]$},
        extra x tick style={tick label style={yshift=-2.5ex}},
        xlabel={Cores\\{}[DoFs $\times 10^6$]},
        xlabel style={align=center, style={yshift=-2ex}},
        ylabel near ticks,
        ylabel style={align=center, text width=4cm},
        ylabel={Time to solution over all continuation steps [min]},
        ]
        \addplot+ table[x=Cores,y=Time] {\bfstable};
      \end{semilogxaxis}
    \end{tikzpicture}
    \caption{Weak scaling of the three-dimensional backwards-facing step.}
    \label{fig:bfs3d-scaling}
  \end{subfigure}
  \caption{Weak scaling of time to solution over all continuation
    steps for both 3D problems.}
\end{figure}
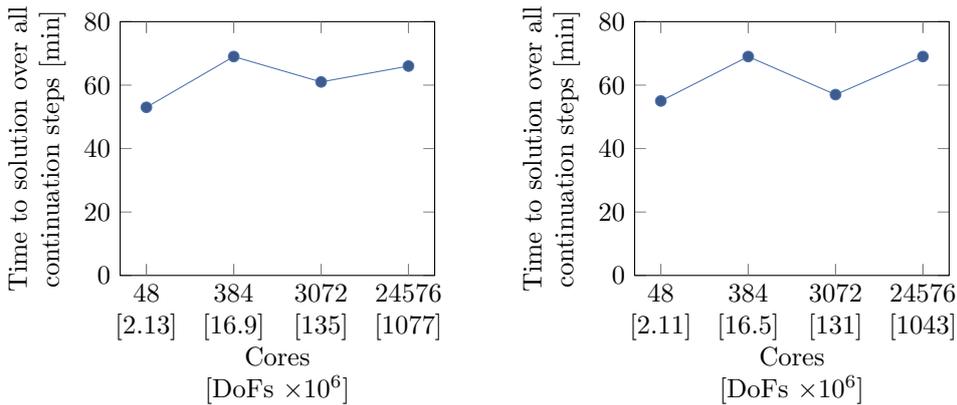

Having seen that the algorithmic scalability of the solver is good,
with well-controlled iteration counts, we now consider the
computational performance. In Figures~\ref{fig:ldc3d-scaling} and
\ref{fig:bfs3d-scaling} we show the weak scaling\footnote{A weak
  scaling test is where the
  number of degrees of freedom per MPI process is held constant while increasing the
  number of processes. Perfect scaling corresponds to a constant time
  to solution as the problem size is increased.} of the total time to
solution over all continuation steps. Both problems show excellent
scalability from 48 to 24576 MPI processes, with the lid-driven cavity
achieving a scaling efficiency of 80\% and the backwards-facing step
79\%.
We attribute the lack of perfect scalability primarily to load imbalance in our mesh distribution. In both
problems, although the mesh has a well-balanced partition of cells,
for the patch smoother to have perfect load balance the number of
vertices owned by each process must also be equal. The partitioning
scheme used does not take this constraint into account, and we observe
that the number of patches per process varies by a factor of 4 over
the partition for the largest problems. The scaling and computational
performance of the code will be improved in future work.

\section{Conclusions and outlook} \label{sec:conclusions}

In this paper we have extended the multigrid method of Benzi, Olshanskii and
Sch\"oberl for the augmented momentum solve arising in the augmented Lagrangian
preconditioner to three dimensions. The prolongation operator proposed by
Sch\"oberl works for the $[\Pthree]^3\mathrm{-}\Pzero$ and
$[\PtwoFB]^3\mathrm{-}\Pzero$ discretizations, while a modification is required
to use the cheaper $[\PoneFB]^3\mathrm{-}\Pzero$ element. We have developed a
new patchwise preconditioner in PETSc and implemented the resulting scheme in
Firedrake. We have demonstrated iteration counts that grow very slowly with
respect to the Reynolds number in both 2D and 3D for problems of up to a billion
degrees of freedom. The code is freely available as open source.

However, this multigrid method is currently tightly coupled to the use of piecewise
constant elements for the pressure for full robustness, and the discretizations considered here do
not represent the divergence-free constraint exactly, which is highly desirable
\cite{john2017}. The key next step is to develop a Reynolds-robust
preconditioner for these discretizations, such as the Scott--Vogelius element
\cite{scott1985}, the Guzm\'an--Neilan modification of Bernardi--Raugel
\cite{guzman2017}, or a $H(\mathrm{div})$-conforming element
\cite{cockburn2006}. It may also be possible to use this solver as a
preconditioner for other discretisations, in a similar manner to the
modified multigrid schemes studied in \cite{john2001}.

\section*{Code availability}
For reproducibility, we cite archives of the exact software versions
used to produce the results in this paper. All major Firedrake
components have been archived on Zenodo~\citep{zenodo/Firedrake-20190617.0}. An
installation of Firedrake with components matching those used to
produce the results in this paper can by obtained following the
instructions at \url{https://www.firedrakeproject.org/download.html}
with\vspace{0.2cm}

{\small
\begin{verbatim}
export PETSC_CONFIGURE_OPTIONS="--download-superlu --download-superlu_dist \
   --with-cxx-dialect=C++11"
python3 firedrake-install --doi 10.5281/zenodo.3247427
\end{verbatim}
}\vspace{0.2cm}

\noindent
The additive Schwarz preconditioner has been incorporated into PETSc
as of version 3.10. The Navier--Stokes solver,
and example files, are available at
\url{https://bitbucket.org/pefarrell/fmwns/}, the version used in the
paper is archived as part of \cite{zenodo/Firedrake-20190617.0}.

\bibliographystyle{siamplain}
\bibliography{literature,zenodo}

\end{document}